\documentclass[11pt]{amsart}
\usepackage{latexsym,amssymb,amsmath}
\input epsf
\newtheorem{proposition}{Proposition}[section]

\newtheorem{corollary}{Corollary}[section]

\def\tdim{\text{dim}\,}

\def\tbase{\text{Base}\,}

\def\BC{\mathbb C}\def\BF{\mathbb F}\def\BO{\mathbb O}
\def\BA{\mathbb A}
\def\BP{\mathbb P}

\def\fr{\mathfrak r}\def\fz{\mathfrak z}

\def\fgl{\mathfrak g\mathfrak l}
\def\ta{\tilde{\alpha}}
\def\hd{, \hdots ,}

\def\cO{{\cal O}}
\def\CC{\mathbb C}
\def\RR{\mathbb R}
\def\HH{\mathbb H}
\def\AA{{\mathbb A}}
\def\BB{{\mathbb B}}
\def\OO{\mathbb O}

\def\ZZ{\mathbb Z}

\def\11{\mathbf 1}
\def\PP{\mathbb P}

\def\FF{\mathbb F}

\def\fh{{\mathfrak h}}
\def\fr{{\mathfrak r}}

\def\fd{{\mathfrak d}}
\def\fsl{{\mathfrak {sl}}}
\def\fsp{{\mathfrak {sp}}}
\def\fspin{{\mathfrak {spin}}}
\def\fso{{\mathfrak {so}}}
\def\fco{{\mathfrak {co}}}
\def\fe{{\mathfrak e}}
\def\ff{{\mathfrak f}}

\def\fg{{\mathfrak g}}
\def\tfg{{\tilde{\mathfrak g}}}
\def\fm{{\mathfrak m}}

\def\fp{{\mathfrak p}}

\def\fk{{\mathfrak k}}
\def\ft{{\mathfrak t}}
\def\fl{{\mathfrak l}}
\def\etimes{\hspace{-1mm}\times\hspace{-1mm}}

\def\a{\alpha}
\def\ta{\tilde{\alpha}}
\def\o{\omega}

\def\b{\beta}

\def\s{\sigma}

\def\th{\theta}

\def\ot{{\mathord{\,\otimes }\,}}
\def\op{{\mathord{\,\oplus }\,}}

\def\der{\mathcal Der}

\def\lra{{\mathord{\;\longrightarrow\;}}}
\def\ra{{\mathord{\;\rightarrow\;}}}

\def\tr{{\rm trace}\;}
\def\dim{{\rm dim}\;}

\def\cJ{{\mathcal J}}

\def\cO{{\mathcal O}}
\def\CC{\mathbb C}
\def\RR{\mathbb R}
\def\HH{\mathbb H}
\def\AA{{\mathbb A}}
\def\BB{{\mathbb B}}
\def\OO{\mathbb O}

\def\ZZ{\mathbb Z}

\def\11{\mathbf 1}
\def\PP{\mathbb P}

\def\FF{\mathbb F}

\def\fh{{\mathfrak h}}

\def\fsl{{\mathfrak {sl}}}
\def\fsp{{\mathfrak {sp}}}
\def\fspin{{\mathfrak {spin}}}
\def\fso{{\mathfrak {so}}}
\def\fe{{\mathfrak e}}

\def\ff{{\mathfrak f}}

\def\fz{{\mathfrak z}}

\def\fg{{\mathfrak g}}

\def\fp{{\mathfrak p}}
\def\fk{{\mathfrak k}}
\def\ft{{\mathfrak t}}
\def\fl{{\mathfrak l}}

\def\a{\alpha}
\def\ta{\tilde{\alpha}}
\def\o{\omega}

\def\b{\beta}

\def\s{\sigma}

\def\th{\theta}

\def\ot{{\mathord{\,\otimes }\,}}
\def\op{{\mathord{\,\oplus }\,}}

\def\lra{{\mathord{\;\longrightarrow\;}}}
\def\ra{{\mathord{\;\rightarrow\;}}}

\def\tr{{\rm trace}\;}
\def\dim{{\rm dim}\;}

\newcommand\rem{{\medskip\noindent {\em Remark}.}\hspace{2mm}}

\begin{document}
\title{Series of Nilpotent orbits}
\author{J.M. Landsberg, Laurent Manivel, Bruce W. Westbury}
\date{September  2003}

\begin{abstract}
We organize the nilpotent orbits in the exceptional complex Lie algebras
into series 
and show that within each series the dimension of the
orbit is a linear function of the natural parameter 
$a=1, 2, 4, 8$, respectively for $\ff_4,\fe_6,\fe_7,\fe_8$.
We observe
similar regularities for the centralizers of nilpotent 
elements in a series and  graded  components in the
associated grading of the ambient Lie
algebra. More strikingly, we observe that for $a\geq 2$ 
the numbers of $\FF_q$-rational 
points on the nilpotent orbits of a given series are 
given by polynomials which
have   uniform expressions in terms of $a$.
This even remains true for the
degrees of the unipotent characters associated 
to these series through the Springer correspondence.
We make similar observations for the series arising from the other
rows of Freudenthal's magic chart and some observations about
the general organization of nilpotent orbits, including
the description of and
dimension formulas for several universal nilpotent orbits (universal in the
sense that they occur in almost all simple Lie algebras).
\end{abstract}

\maketitle

\section{Introduction}

\subsection{Main results}
In this paper 
we explore consequences of the Tits-Freudenthal construction and its variant, the triality model,
for nilpotent orbits in  the exceptional complex simple Lie algebras. 
Both models produce a Lie algebra $\fg(\AA,\BB)$ 
from a pair $\AA, \BB$ of real normed algebras. When $\BB=\OO$, one obtains the exceptional
Lie algebras $\ff_4, \fe_6, \fe_7, \fe_8$, parametrized by the dimension $a=1, 2, 4, 8$ of 
$\AA$. In \cite{LMtrial}, the first two authors used the triality model to explain rather mysterious 
formulas obtained by Deligne for the dimensions of 
certain series of representations of the exceptional Lie
algebras. In this paper, we show that the use of the parameter $a$ leads to
several interesting observations for nilpotent orbits in the exceptional Lie algebras. 

\smallskip
Let us begin with any nilpotent orbit $\cO$ in $\ff_4$. Thanks to the natural embeddings
$\fso_8\subset\ff_4\subset\fe_6\subset\fe_7\subset\fe_8$, we obtain a series of orbits $\cO_a$ in these
Lie algebras. Their weighted Dynkin diagrams can be obtained in the following way: 
each series of orbits is defined by a weight of $\fso_8$, which, through the triality
model,  defines a weight, thus
a weighted Dynkin diagram, for each exceptional Lie algebra. This is  how   
we proved and generalized the dimension formulas of Deligne for 
series of representations whose highest weights came from $\fso_8$
in \cite{LMtrial}.
We prove, or check from the tables compiled in \cite{carter}, that 
\begin{itemize}
\item the dimension of $\cO_a$ is a linear function of $a$, 
\item the stabilizers of points in $\cO_a$ have unipotent radicals of dimension again linear
in $a$, while their reductive parts organize into simple series, 
\item the closure of $\cO_a$ can be desingularized by a homogeneous vector bundle, 
whose dimensions of the base and of the fiber, are both linear in $a$, 
\item the number of $\FF_q$-rational points on $\cO_a$, for large $q$, is given by 
a polynomial in $q$ with a uniform expression in $a$, 
\item the unipotent characters of the finite groups of exceptional Lie type, associated
to the orbits $\cO_a$  through the Springer correspondence, have degrees given by 
polynomials that, when suitably expressed as rational functions, have uniform 
expressions in $a$. 
\end{itemize}
This last fact, which is true only for $a\ge 2$, is the most mysterious 
observation of this paper, and we would like very much 
to have a theoretical  explanation. 
   
\smallskip
We observe similar phenomena for the other lines of Freudenthal's square i.e., for the
series of Lie algebras $\fg(\AA,\BB)$ when $\BB$ is $\RR, \CC$ or $\HH$, and also for 
the classical Lie algebras. In fact, a few orbits are universal (or almost universal), 
in the sense that they appear in every (or almost every) simple Lie algebra. We discuss 
certain properties of these orbits in connection with the work of Vogel and Deligne around
the ``universal Lie algebra'', and also with the more geometric investigations
of \cite{LMmagic}.

\subsection{The Freudenthal-Tits construction}

We recall Tits' construction of the exceptional Lie algebras in terms of 
real normed algebras \cite{tits0}. 

Let $\AA$ a be a real normed algebra, 
so that   $\AA=\RR, \CC, \HH$ or $\OO$, the Cayley algebra, and let $a:={\rm dim}
\AA = 1, 2, 4$ or $8$. The conjugation (i.e., the orthogonal symmetry with respect 
to the unit element) will be denoted by $u\mapsto u^*$. The subspace of $\AA$ defined by the
equation $u^*=-u$ is the orthogonal $Im\AA$ of the unit element. Let $\cJ_3(\AA)$ denote the 
Jordan algebra of Hermitian matrices of order three, with coefficients in $\AA$. 
The subspace of traceless matrices is denoted $\cJ_3(\AA)_0$. 
 
Now let $\AA$ and $\BB$ be two real normed algebras, and let 
$$\fg(\AA,\BB)=Der\AA\times Der \cJ_3(\BB)\op (Im\AA \ot  \cJ_3(\BB)_0).$$ 
There is a natural structure of $\ZZ_2$-graded Lie algebra on $\fg(\AA,\BB)$. 

A useful variant of this construction is the triality model, first discovered by
Allison \cite{allison}, and recently rediscovered by several authors
(see e.g., \cite{LMtrial}). 
Define the {\sl triality algebra}
$$\ft(\AA)=\{\th=(\th_1,\th_2,\th_3)\in \fso(\AA)^3,\;\;
\th_3(xy)=\th_1(x)y+x\th_2(y)\;\;\;\forall x,y\in\AA\}.$$
We have $\ft(\RR)=0$,  $\ft(\CC)=\RR^2$, $\ft(\HH)=\fso_3\etimes\fso_3\etimes\fso_3$
and $\ft(\OO)=\fso_8$. 

For $\AA$ and $\BB$ two real normed algebras, let 
$$\tfg(\AA,\BB) = \ft(\AA)\times\ft(\BB)\oplus (\AA_1\ot\BB_1)
\oplus (\AA_2\ot\BB_2) \oplus (\AA_3\ot\BB_3). $$
Then there is a natural structure of $\ZZ_2\times\ZZ_2$-graded semi-simple 
Lie algebra on $\tfg(\AA,\BB)$. 

In what follows we will work over the complex numbers and complexify the whole construction 
without changing notations. We just have a new conjugation map $x\ra\overline{x}$ in $\OO$
(which now denotes the complexified Cayley algebra) such that $\overline{xy}=
\overline{x}\times\overline{y}$. 

The result of both constructions is Freudenthal's magic square:

\medskip
\begin{center}\begin{tabular}{|c|cccc|} \hline 
 & $\RR$ & $\CC$ & $\HH$ & $\OO$ \\ \hline
$\RR$ &  $\fsl_2$ & $\fsl_3$ & $\fsp_6$ & $\ff_4$ \\  
$\CC$ & $\fsl_3$ & $\fsl_3\etimes \fsl_3$ & $\fsl_6$ & $\fe_6$ \\ 
$\HH$ & $\fsp_6$ & $\fsl_6$ & $\fso_{12}$ & $\fe_7$  \\
$\OO$ & $\ff_4$ & $\fe_6$ & $\fe_7$ & $\fe_8$ \\ \hline
\end{tabular}\end{center}
\medskip

\subsection{The exceptional series}

Letting $\BB=\OO$ in the magic square, we obtain exceptional Lie algebras of types 
$\ff_4, \fe_6, \fe_7$ and $\fe_8$. 

\smallskip
An important fact for what follows is the observation (\cite{LMtrial}, page 68) 
that there is a preferred cone 
${\mathcal C}$ in the weight lattice of $\ft(\OO)=\fso_8$ defined by the condition that a 
weight in ${\mathcal C}$ is dominant and integral when considered as a weight of 
each of the four Lie algebras $\fg(\AA,\OO)\supset\fso_8$. This cone is generated
by the four following weights of $\fso_8$:
$$\o(\fg)=\o_2,\quad \o(\fg_2)=\o_1+\o_3+\o_4, \quad \o(\fg_3)=2\o_1+2\o_3, \quad 
\o(\fg_Q)=2\o_1.$$

\smallskip\noindent (The representations we denote by $\fg, \fg_2, \fg_3, \fg_Q$ are
denoted $X_1, X_2, X_3, Y_2^*$ in \cite{del,LMtrial}.)

The following table contains the expressions of these four weights in terms of the 
fundamental weights of each exceptional Lie algebra. 

$$\begin{array}{rcccc}
 & \ff_4 & \fe_6 & \fe_7 & \fe_8 \\
\o(\fg) & 
\setlength{\unitlength}{2mm}
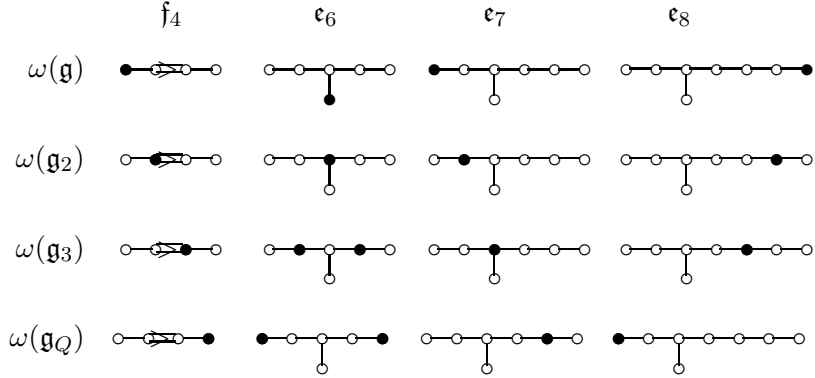
\begin{picture}(7,3)(0,0)
\multiput(0,0)(2,0){4}{$\circ$}
\put(0,0){$\bullet$}
\multiput(0.8,.5)(4,0){2}{\line(1,0){1.4}} 
\multiput(2.55,.3)(0,.5){2}{\line(1,0){1.7}} 
\put(2.5,0){$>$}
\end{picture} & 
\setlength{\unitlength}{2mm}
\begin{picture}(8,3)(0.2,0)
\multiput(0,0)(2,0){5}{$\circ$}
\multiput(0.7,.5)(2,0){4}{\line(1,0){1.4}} 
\put(4,-2){$\bullet$}
\put(4.45,-1.3){\line(0,1){1.4}}
\end{picture} & 
\setlength{\unitlength}{2mm}
\begin{picture}(9,3)(0,0)
\multiput(2,0)(2,0){5}{$\circ$}
\put(0,0){$\bullet$}
\multiput(0.7,.5)(2,0){5}{\line(1,0){1.4}} 
\put(4,-2){$\circ$}
\put(4.45,-1.25){\line(0,1){1.4}}
\end{picture} & 
\setlength{\unitlength}{2mm}
\begin{picture}(10,3)(-1,0)
\multiput(0,0)(2,0){6}{$\circ$}
\put(12,0){$\bullet$}
\multiput(0.7,.55)(2,0){6}{\line(1,0){1.4}} 
\put(4,-2){$\circ$}
\put(4.45,-1.25){\line(0,1){1.4}}
\end{picture} \\
 & & & & \\
\o(\fg_2) &
\setlength{\unitlength}{2mm}
\begin{picture}(7,3)(0,0)
\multiput(0,0)(2,0){4}{$\circ$}
\put(2,0){$\bullet$}
\multiput(0.8,.55)(4,0){2}{\line(1,0){1.4}} 
\multiput(2.55,.3)(0,.5){2}{\line(1,0){1.7}} 
\put(2.5,0){$>$}
\end{picture}
& 
\setlength{\unitlength}{2mm}
\begin{picture}(8,3)(0.2,0)
\multiput(0,0)(2,0){5}{$\circ$}
\put(4,0){$\bullet$}
\multiput(0.7,.55)(2,0){4}{\line(1,0){1.4}} 
\put(4,-2){$\circ$}
\put(4.45,-1.25){\line(0,1){1.4}}
\end{picture}
& 
\setlength{\unitlength}{2mm}
\begin{picture}(9,3)(0,0)
\multiput(0,0)(2,0){6}{$\circ$}
\multiput(0.7,.55)(2,0){5}{\line(1,0){1.4}} 
\put(4,-2){$\circ$}
\put(2,0){$\bullet$}
\put(4.45,-1.25){\line(0,1){1.4}}
\end{picture}
& 
\setlength{\unitlength}{2mm}
\begin{picture}(10,3)(-1,0)
\multiput(0,0)(2,0){7}{$\circ$}
\multiput(0.7,.55)(2,0){6}{\line(1,0){1.4}} 
\put(4,-2){$\circ$}
\put(10,0){$\bullet$}
\put(4.45,-1.25){\line(0,1){1.4}}
\end{picture} \\
   & & & \\

\o(\fg_3) &
\setlength{\unitlength}{2mm}
\begin{picture}(7,3)(0,0)
\put(4,0){$\bullet$}
\multiput(0,0)(2,0){4}{$\circ$}
\multiput(0.8,.5)(4,0){2}{\line(1,0){1.4}} 
\multiput(2.55,.3)(0,.5){2}{\line(1,0){1.7}} 
\put(2.5,0){$>$}
\end{picture}
& 
\setlength{\unitlength}{2mm}
\begin{picture}(8,3)(0.2,0)
\multiput(0,0)(2,0){5}{$\circ$}
\multiput(0.7,.5)(2,0){4}{\line(1,0){1.4}} 
\put(2,0){$\bullet$}\put(6,0){$\bullet$}
\put(4,-2){$\circ$}
\put(4.45,-1.3){\line(0,1){1.4}}
\end{picture}
& 
\setlength{\unitlength}{2mm}
\begin{picture}(9,3)(0,0)
\multiput(0,0)(2,0){6}{$\circ$}
\put(4,0){$\bullet$}
\multiput(0.7,.5)(2,0){5}{\line(1,0){1.4}} 
\put(4,-2){$\circ$}
\put(4.45,-1.3){\line(0,1){1.4}}
\end{picture}
& 
\setlength{\unitlength}{2mm}
\begin{picture}(10,3)(-1,0)
\multiput(0,0)(2,0){7}{$\circ$}
\put(8,0){$\bullet$}
\multiput(0.7,.5)(2,0){6}{\line(1,0){1.4}} 
\put(4,-2){$\circ$}
\put(4.45,-1.3){\line(0,1){1.4}}
\end{picture} \\
 &  & & & \\
\o(\fg_Q) &
\setlength{\unitlength}{2mm}
\begin{picture}(8,3)(0,0)
\multiput(0,0)(2,0){4}{$\circ$}
\put(6,0){$\bullet$}
\multiput(0.8,.55)(4,0){2}{\line(1,0){1.4}} 
\multiput(2.55,.3)(0,.5){2}{\line(1,0){1.7}} 
\put(2.5,0){$>$}
\end{picture}
& 
\setlength{\unitlength}{2mm}
\begin{picture}(9,3)(0.2,0)
\multiput(0,0)(2,0){5}{$\circ$}
\put(0,0){$\bullet$}
\put(8,0){$\bullet$}
\multiput(0.7,.55)(2,0){4}{\line(1,0){1.4}} 
\put(4,-2){$\circ$}
\put(4.45,-1.25){\line(0,1){1.4}}
\end{picture}
& 
\setlength{\unitlength}{2mm}
\begin{picture}(10,3)(0,0)
\multiput(0,0)(2,0){6}{$\circ$}
\multiput(0.7,.55)(2,0){5}{\line(1,0){1.4}} 
\put(4,-2){$\circ$}
\put(8,0){$\bullet$}
\put(4.45,-1.25){\line(0,1){1.4}}
\end{picture}
& 
\setlength{\unitlength}{2mm}
\begin{picture}(11,3)(-1,0)
\multiput(0,0)(2,0){7}{$\circ$}
\put(0,0){$\bullet$}
\multiput(0.7,.55)(2,0){6}{\line(1,0){1.4}} 
\put(4,-2){$\circ$}
\put(4.47,-1.25){\line(0,1){1.4}}
\end{picture}

\end{array}$$

\bigskip

\bigskip

\section{Nilpotent orbits in the exceptional series}

\subsection{Series of nilpotent orbits}
Since $Der\cJ_3(\OO)=\ff_4$ is a subalgebra of $\fg(\AA,\OO)$ for all $\AA$, every nilpotent
orbit in $\ff_4$ defines a nilpotent orbit $\cO_a$ in $\fg(\AA,\OO)$, of the corresponding 
adjoint Lie group. More generally, any element of $\ff_4$ defines a series of orbits in the Lie 
algebras $\fg(\AA,\OO)$.

\begin{proposition}
For any element of $\ff_4$, the dimension of its orbit $\cO_a$ in $\fg(\AA,\OO)$, 
is a linear function of $a$. 
\end{proposition}

\proof Let $X$ belong to $\ff_4$, and let's denote its centralizer by $c(X)\subset\ff_4$. 
The centralizer $c(X)_a$ of $X$ in $\fg(\AA,\OO)$ is 
$Der\AA\times c(X)\op Im\AA\ot k(X)$, where $k(X)\subset
\cJ_3(\OO)_0$ denotes the subspace annihilated by $X$. The codimension of this centralizer is 
obviously a linear function of $a$. Since it is equal to the dimension of the orbit $\cO_a$ 
of $X$ in $\fg(\AA,\OO)$, our claim is proved. \qed

\smallskip
Now we suppose that $X\in\ff_4$ is nilpotent, and we
complete it into a $\fsl_2$-triple $(X,Y,H)$ of $\ff_4$. 
The reductive part of $c(X)_a$ is the centralizer $\fh(a):=c(X,Y,H)_a$ of the full 
$\fsl_2$-triple (\cite{carter}, Proposition 5.5.9). Moreover, the decomposition of the 
adjoint action of $H$ into eigenspaces is 
$$\fg(\AA,\OO)=\bigoplus_{i\in\ZZ}\fg(a,i),$$
with $[\fg(a,i),\fg(a,j)]\subset \fg(a,i+j)$. In particular, $\fg(a,0)$ is a subalgebra, 
and each $\fg(a,i)$ is a $\fg(a,0)$-module. Note that  $\fg(a,0)$ contains $\fh(a)$.

\begin{proposition}
For every nilpotent orbit $\cO_1$ in $\ff_4$, let again $\cO_a$ denote the corresponding 
series of nilpotent orbits in $\fg(\AA,\OO)$. 
The dimension of the nilpotent radical $\fr(a)$ of the stabilizer of an element of $\cO_a$ 
is a linear function of $a$. 
For any  $i\ne 0$, the dimension of the $i$-th part $\fg(a,i)$ 
of the induced gradation of $\fg(\AA,\OO)$ 
is a linear function of $a$. 
\end{proposition}

\proof Let $X\in\cO_1\subset\ff_4$ be nilpotent, and  
$(X,Y,H)$ a $\fsl_2$-triple  of $\ff_4$.  If 
$k(X,Y,H)=k(X)\cap k(Y)\cap k(H)$, the centralizer 
of the $\fsl_2$-triple is 
$$c(X,Y,H)_a=Der\AA\times c(X,Y,H)_1\op Im\AA\ot k(X,Y,H),$$
whose codimension in $c(X)_a$ is a linear function of $a$.
Since this is the reductive part $\fh(a)$ of this centralizer,  
its codimension is equal to the 
dimension of the nilpotent radical $\fr(a)$ of $c(X)_a$, 
and our first claim is proved. 

For the second claim, we just note that for $i\ne 0$, $\fg(a,i)=\fg(0,i)\op Im\AA\otimes\fk(i)$, 
where $\fk(i)\subset \cJ_3(\OO)_0$ is the $i$-th eigenspace of the $H$-action. 
(We thank E. Vinberg for these observations.) \qed

\medskip
\rem Note that since $\fh(a)$ centralizes the $\fsl_2$-triple, $\fh(a)\times\fsl_2$ is naturally 
a subalgebra of $\fg(\AA,\OO)$ which can be decomposed into
$$\fg(\AA,\OO)=\bigoplus_{k\ge 0}\fg^*(a,k)\ot [k],$$
where $[k]$ denotes the irreducible $\fsl_2$-module of dimension $k+1$, and 
$\fg^*(a,k)$ is a $\fh(a)$-module. In particular, $\fg^*(a,0)=\fh(a)$. By elementary
properties of the representation theory of $\fsl_2$, the dimension of $\fg^*(a,k)$
is $\dim \fg(a,k)-\dim \fg(a,k+2)$, and is again a linear function of $a$ for $k\ne 0$. 

\medskip
Recall that the nilpotent orbits can be classified by combinatorial data as follows: 
If $X$ belongs to some nilpotent orbit $\cO$, we include it into a $\fsl_2$-triple $(X,Y,H)$.
The semi-simple element $H$ can be supposed to belong to a given Cartan subalgebra $\ft$, and a set 
of simple roots $\Delta$ can be chosen such that $\a(H)$ is a non negative integer for all
$\a\in\Delta$. The collection of these integers, or the corresponding weighted Dynkin diagram, 
uniquely defines the nilpotent orbit $\cO$. 

\smallskip To understand the weighted Dynkin diagrams of a series $\cO_a$ of nilpotent orbits
in the exceptional Lie algebras, it is convenient to use the triality model $\tfg(\AA,\OO)$ 
rather than the more classical Tits-Freudenthal construction. Beginning with $\AA=\RR$, 
we have 
$$\ff_4=\tfg(\RR,\OO)=\fso_8\op\OO_1\op\OO_2\op\OO_3.$$
A Cartan subalgebra $\ft$ of $\ff_4$ can be chosen inside $\fso_8$. We use the notations
of \cite{bou} for the root system of $\fso_8$ and choose the same simple roots. The roots of 
$\ff_4$ are then given by those of $\fso_8$, plus the weights of the three 
inequivalent eight-dimensional 
representations $\OO_1$, $\OO_2$, $\OO_3$. We get a set of positive roots by choosing a linear 
form on $\ft^*$ of the form $\ell=\ell_1\a_1^*+\ell_2\a_2^*+\ell_3\a_3^*+\ell_4\a_4^*$, with 
$\ell_1>\ell_2>\ell_3>\ell_4>0$. The three representations  $\OO_1$, $\OO_2$, $\OO_3$ have 
highest weight $\o_1$, $\o_3$, $\o_4$ respectively, and their minimal weights on which 
$\ell$ is positive are $\phi_1=\o_3-\o_4$, $\phi_2=\o_1-\o_4$, $\phi_3=\o_1-\o_3$ respectively. 
The simple roots of $\ff_4$ must either be simple roots of $\fso_8$, or among these
three minimal weights. Since $\phi_3=\phi_1+\phi_2$, $\a_3=\a_4+2\phi_1$ and $\a_1=\a_3+2\phi_4$, 
the simple roots of $\ff_4$ must be  $\a_2, \a_4,\phi_1$ and $\phi_2$. Note that our four preferred 
weights $\o(\fg)$,  $\o(\fg_2)$,  $\o(\fg_3)$,  $\o(\fg_Q)$ of $\fso_8$ provide us with the dual basis.

\smallskip
Now let $\AA$ be any real normed algebra (complexified). A Cartan subalgebra 
of $\fg(\AA,\OO)$ is given by the sum of the Cartan subalgebra $\ft$ of $\ft(\OO)=\fso_8$, 
and a Cartan subalgebra of $\ft(\AA)$. The root system of $\fg(\AA,\OO)$ is the union of 
the roots systems of $\fso_8$ and $\ft(\AA)$, plus the weights of the form $\mu+\nu$, 
for $\mu$ a weight of some $\OO_i$, and $\nu$ a weight of $\AA_i$. The positive roots can 
be chosen to be the positive roots of $\fso_8$ and $\ft(\AA)$, plus the weights $\mu+\nu$
for which $\ell (\mu)>0$. The simple roots of $\fg(\AA,\OO)$ are then either simple roots
of $\fso_8$, of $\ft(\AA)$ (we denote them by $\a'_j$), or some of the $\phi_i-\o'_i$, 
where $\o'_i$ is the highest weight of $\AA_i$ (and $-\o'_i$ its lowest weight, since 
$\AA_i$ is self-dual). Since $S^2\AA_i$ contains the trivial representation, $2\o'_i$ 
must belong to the root lattice of $\ft(\AA)$, as well as $\o'_1+\o'_2+\o'_3$ because 
there is an equivariant map $\AA_1\ot\AA_2\ra\AA_3$. We easily deduce that exactly 
as in the case of $\ff_4$, $\phi_3-\o'_3$, $\a_1$ and $\a_3$ cannot be simple roots. 
The simple roots of $\tfg(\AA,\OO)$ are therefore given by $\a_2, \a_4$, the $\a'_j$'s, 
$\phi_1-\o'_1$ and $\phi_2-\o'_2$.

\smallskip
For a $\fsl_2$-triple $(X,Y,H)$ in $\ff_4=\fg(\RR,\OO)$,
defining a nilpotent orbit $\cO_a$ in $\tfg(\AA,\OO)$,  the labels of the corresponding 
Dynkin diagram will be  $\a_2(H)$,  $\a_4(H)$, $\a'_j(H)=0$, $\phi_1(H)$ and $\phi_2(H)$, 
i.e., exactly the same labels as those of $\cO_1$, plus some zeros on the simple roots coming 
from $\ft(\AA)$.  We conclude:

\begin{proposition}
Let the nilpotent orbit $\cO_1$ in $\ff_4$ define a series $\cO_a$ of nilpotent orbits 
in the exceptional Lie algebras. Suppose that the weighted Dynkin diagram of $\cO_1$ defines the
weight $p\o(\fg)+q\o(\fg_2)+r\o(\fg_3)+s\o(\fg_Q)$. Then this remains true for the weighted Dynkin 
diagrams of each of the nilpotent orbits $\cO_a$. 
\end{proposition}
 
\pagebreak
We encode the corresponding series by the symbol $\fg^p\fg_2^q\fg_3^r\fg_Q^s$. With this convention, 
the Hasse diagram of nilpotent orbits in $\ff_4$ (see e.g., \cite{carter}, p. 440), 
is given by the following picture: 

{\small
\setlength{\unitlength}{2.5mm}
\begin{picture}(20,55)(-16,-1)

\put(10,49){$\fg^2\fg_2^2\fg_3^2\fg_Q^2$}
\put(11,45){$\fg^2\fg_2^2\fg_Q^2$}
\put(12.5,46.3){\line(0,1){2}}
\put(11,41){$\fg_2^2\fg_Q^2$}
\put(12.5,42.3){\line(0,1){2}}
\put(5,37){$\fg^2\fg_2^2$}
\put(16,37){$\fg\fg_3\fg_Q^2$}
\put(11.8,33){$\fg_2^2$}
\put(7.8,38.5){\line(1,1){1.5}}
\put(14.3,40){\line(1,-1){1.5}}
\put(11.5,29){$\fg\fg_3$}
\put(9.8,34.5){\line(-1,1){1.5}}
\put(13.8,34.5){\line(1,1){1.5}}
\put(4.5,25){$\fg_2\fg_Q$}
\put(12.5,30.3){\line(0,1){2}}
\put(16.5,25){$\fg^2\fg_Q$}
\put(7.8,26.5){\line(1,1){1.5}}
\put(14.3,28){\line(1,-1){1.5}}
\put(17.5,21){$\fg_3$}
\put(8,24){\line(3,-1){7.5}}
\put(5,17){$\fg_Q^2$}
\put(5.5,18.3){\line(0,1){6}}
\put(18,18.3){\line(0,1){2}}
\put(18,22.3){\line(0,1){2}}
\put(17.5,17){$\fg^2$}
\put(9.8,14.5){\line(-1,1){1.5}}
\put(13.8,14.5){\line(1,1){1.5}}
\put(11.6,13){$\fg_2$}
\put(12.5,10.2){\line(0,1){2}}
\put(11.5,9){$\fg_Q$}
\put(12.5,6.2){\line(0,1){2}}
\put(12,5){$\fg$}
\put(12.5,2.2){\line(0,1){2}}
\put(12.2,.6){$0$}

\end{picture} }

\centerline{{\small\sl Hasse diagram of nilpotent orbits in $\ff_4$}}

\bigskip

\smallskip\noindent {\sl Example}. The series of nilpotent orbits $\fg\fg_3\fg_Q^2$
will be given by the following four weighted Dynkin diagrams:

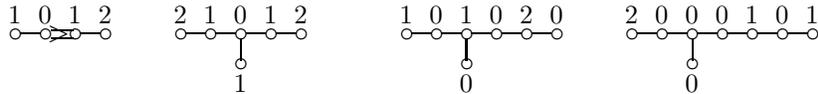
\setlength{\unitlength}{2mm}
\begin{picture}(48,7)(-22,-2)

\multiput(-11,0)(2,0){4}{$\circ$}
\multiput(-10.2,.5)(4,0){2}{\line(1,0){1.3}} 
\multiput(-8.1,.2)(0,.5){2}{\line(1,0){1.4}} 
\put(-8.4,0){$>$}
\put(-11,1.2){{\small $1$}}
\put(-9,1.2){{\small $0$}}
\put(-7,1.2){{\small $1$}}
\put(-5,1.2){{\small $2$}}

\multiput(0,0)(2,0){5}{$\circ$}
\multiput(0.8,.5)(2,0){4}{\line(1,0){1.3}} 
\put(4,-2){$\circ$}
\put(4.45,-1.3){\line(0,1){1.4}}
\put(0,1.2){{\small $2$}}
\put(2,1.2){{\small $1$}}
\put(4,1.2){{\small $0$}}
\put(6,1.2){{\small $1$}}
\put(8,1.2){{\small $2$}}
\put(4,-3.4){{\small $1$}}

\multiput(15,0)(2,0){6}{$\circ$}
\multiput(15.8,.5)(2,0){5}{\line(1,0){1.3}} 
\put(19,-2){$\circ$}
\put(19.45,-1.3){\line(0,1){1.4}}
\put(15,1.2){{\small $1$}}
\put(17,1.2){{\small $0$}}
\put(19,1.2){{\small $1$}}
\put(21,1.2){{\small $0$}}
\put(23,1.2){{\small $2$}}
\put(25,1.2){{\small $0$}}
\put(19,-3.4){{\small $0$}}

\multiput(30,0)(2,0){7}{$\circ$}
\multiput(30.8,.5)(2,0){6}{\line(1,0){1.3}} 
\put(34,-2){$\circ$}
\put(34.5,-1.3){\line(0,1){1.4}}
\put(30,1.2){{\small $2$}}
\put(32,1.2){{\small $0$}}
\put(34,1.2){{\small $0$}}
\put(36,1.2){{\small $0$}}
\put(38,1.2){{\small $1$}}
\put(40,1.2){{\small $0$}}
\put(42,1.2){{\small $1$}}
\put(34,-3.4){{\small $0$}}
\end{picture}

\medskip

\subsection{Series of stabilizers}

For each series $\cO_a$, we proved in Proposition 2.2 that the 
codimension of the centralizer, and the dimension of the nilpotent radical $\fr(a)$, 
are linear  functions of $a$. 
In this section we provide explicit data for each series of orbits. 
We also give the reductive parts $\fh(a)$ of these centralizers, and observe they 
organize into series of Lie algebras. Most of these
are either given by the other series $\fg(\AA,\BB)$ of Freudenthal's square,
the derivation algebras $\der\AA$, the triality algebras $\ft(\AA)=
\der\AA\op 2{\rm Im}\AA$, or the intermediate series $\fl(\AA)=\der\AA\op 
{\rm Im}\AA$ of Barton and Sudbery (\cite{bs}, page 13).
 
Another series that appears is the {\it inf-Severi series} $\fk(\AA)$. 
It has two preferred representations $V(a)$ and $W(a)$,
respectively of dimensions $2a$ and $a+2$. Geometrically, let $X(a)$ be one
of the four Severi varieties, which is homogeneous under the action of the 
adjoint group of $\fg(\AA,\CC)$ \cite{LMpop}. Then $\fk(\AA)$ is the reductive part of the 
Lie algebra of the stabilizer of a point in $X(a)$, $V(a)$ is the isotropy representation, 
and $W(a)$ is the complement of the Cartan square of $V(a)^*$ in $S^2V(a)^*$
(except when $a=1$, in which case it is equal to this Cartan square).

These series of Lie algebras are given by:
$$\begin{array}{rcccc}
\AA & \RR & \CC & \HH & \OO \\
\der\AA&  0 & 0 & \fsl_2 & \fg_2 \\
\fl(\AA) & 0 & \CC & 2\fsl_2 & \fspin_7 \\
\ft(\AA) &  0 & 2\CC & 3\fsl_2 & \fspin_8 \\
\fk(\AA) & \fsl_2 &  \fsl_2\times\fgl_2 &  \fsl_2\times\fsl_4 &  \fspin_{10} 
\end{array}$$

\smallskip
Most of the data below have been 
gathered from the tables in \cite{carter}. 
We refer to each series of orbits by its label $\fg^p\fg_2^q\fg_3^r\fg_Q^s$.
Then we provide the series of labels 
used in the tables of \cite{carter}: in general four of them, encoding the four orbits
in $\ff_4$, $\fe_6$,  $\fe_7$,  $\fe_8$; sometimes five, when the series comes from 
$\fso_8\subset\ff_4$, in which case we also provide the partition of $8$ encoding the corresponding 
orbit (actually sometimes a trialitarian triple of orbits) in $\fso_8$, which corresponds 
to $a=0$. 

\smallskip\noindent {\sl Remark}.  If an $\fso_8$ orbit is symmetric about its folding, it
also occurs in $\fg_2$, and its dimension is given by the same formula with $a=-2/3$. 
This occurs for the orbits
labelled $\fg,\fg_2,\fg^2,\fg^2\fg_2^2$. Similarly, the  formulas
for $\fg$ extend to both $\fsl_2$ and $\fsl_3$ with $a=-4/3$ and $a= -1$, respectively, 
and $\fg^2$ extends also to $\fsl_3$. That these Lie algebras should be incorporated in the
exceptional series was already observed in \cite{del}.
\hspace{-5mm}
$$\begin{array}{ll}
\fg &  \dim\cO_a = 6a+10 \\
 & \dim\fr(a)=6a+9 \\
\, [(2^21^4),A_1,A_1,A_1,A_1]\hspace{3cm} &  
\fh(a) = 3\fsl_2, \; \fsp_6,\; \fsl_6, \; \fso_{12}, \; \fe_7\hspace{25mm}
\end{array}$$

\noindent This is the minimal nilpotent orbit, the cone over the adjoint variety. 
Here $\fh(a)=\fg(\AA,\HH)$, $\fg(a,0)=\fg(\AA,\HH)\times\CC$, $\fg(a,1)=\fz_2(\AA)$, 
the Zorn representation (see for example \cite{LMpop}), and $\fg(a,2)=\CC$. 

$$\begin{array}{ll}
\fg_Q & \dim\cO_a = 10a+12 \\
 & \dim\fr(a)=9a+6 \\
\, [(2222),\tilde{A_1},2A_1,2A_1,2A_1] \hspace{2cm} 
& \fh(a) = \fso_5, \; \fsl_4,\; \fco_7, \; \fso_9\times\fsl_2, 
\; \fso_{13} \qquad\qquad
\end{array}$$

\noindent We denoted by $\fco_n=\fso_n\times\CC$ the conformal Lie algebra. 
Here $\fg(a,0) = \fco_3, \; \fco_7,\; \fco_8, \; \fco_{10}\times\fsl_2, \; \fco_{14}$
respectively, 
$\fg(a,1)$ is a spin representation of dimension $8a$, and for $a>0$, $\fg(a,2)$ is the 
standard vector representation, of dimension $a+6$.

$$\begin{array}{ll}
\hspace{-2mm}\fg_2 & \dim\cO_a = 12a+16 \\
  & \dim\fr(a)=9a+9 \\
\hspace{-2mm} [(3221),A_1+\tilde{A_1},3A_1,3A_1,3A_1] \qquad 
 & \fh(a) = \fsl_2, \; 2\fsl_2,\; \fsl_2\times\fsl_3, \; \fsl_2\times\fsp_6, \; \fsl_2\times\ff_4
\end{array}$$

\noindent This is the series of orbits discussed by Panyushev in \cite{pan}. 
Here  $\fh(a)=\fsl_2\times\fg(\AA,\RR)$ and 
$\fg(a,0)=\fg(\AA,\CC)\times\fgl_2$. If $U$ denotes the natural two-dimensional representation 
of this $\fgl_2$, we have $\fg(a,1)=\cJ_3(\AA)\otimes U$, $\fg(a,2)=\cJ_3(\AA)$, 
and $\fg(a,3)=U$. 
 
$$\begin{array}{ll}
\fg^2  & \dim\cO_a = 12a+18 \\
  & \dim\fr(a)=6a+8 \\
 \, {\small [(3311),A_2,A_2,A_2,A_2]}\hspace{2.5cm}
 & \fh(a) = 2\CC, \;  \fsl_3,\; 2\fsl_3, \; \fsl_6, \; \fe_6  \hspace*{3cm}
\end{array}$$
\noindent This is the $a=2$ line of the Freudenthal square, that is $\fh(a)=\fg(\AA,\CC)$.  
Moreover, since this is the orbit $\fg^2$, the induced grading is the same as in the 
case of the minimal nilpotent orbit, with indices doubled: $\fg(a,0)=\fg(\AA,\HH)\times\CC$, 
$\fg(a,1)=0$, $\fg(a,2)=\fz_2(\AA)$, $\fg(a,3)=0$, $\fg(a,4)=\CC$.

$$\begin{array}{ll}
\hspace{-3mm}\fg_3 & \dim\cO_a = 16a+18 \\
 &  \dim\fr(a)=9a+6 \\
\hspace{-2mm} [A_2+\tilde{A_1},A_2+2A_1,A_2+2A_1,A_2+2A_1]\quad  
 & \fh(a) = \fsl_2, \; \fgl_2,\; 3\fsl_2, \; \fsl_2\times\fso_7 \hspace{1cm}
\end{array}$$
\noindent Here $\fh(a)=\fsl_2\times\fl(a)$. Moreover, 
$\fg(a,0)=\fgl_3\times\fk(\AA)$, where $\fk(\AA)$ is the inf-Severi series discussed above. 
Let $U$ denote the natural representation of $\fgl_3$. 
Then $\fg(a,1)=U\otimes V(a)$, 
$\fg(a,2)=U^*\ot W(a)$, $\fg(a,3)=V(a)$ and $\fg(a,4)=U$. 
$$\begin{array}{ll}
\hspace{1cm}\fg^2\fg_Q &  \dim\cO_a = 16a+20 \\
 & \dim\fr(a)=5a+5 \\
\hspace{1cm} [(44),B_2,A_3,A_3,A_3]\hspace{3cm} 
 & \fh(a) = 2\CC, \; 2\fsl_2 ,\; \fco_5, \; \fso_7\times\fsl_2, \; \fso_{11}\qquad \quad
\end{array}$$
\noindent Here $\fg(a,0)=2\fgl_2, \fco_5, \fgl_4\times\CC^2, \fgl_2\times\fco_8, 
\fco_{12}\times\CC$ respectively. Moreover, $\fg(a,1)$ and $\fg(a,3)$ have dimension $4a$, 
$\fg(a,2)$ and $\fg(a,4)$ have dimension $a+4$, and $\fg(a,5)$ is one-dimensional.
For $a=1$ we get representations of dimensions $4$ and $5$, in accordance 
with the exceptional isomorphism $\fso_5\simeq\fsp_4$.    
$$\begin{array}{ll}
\hspace{5mm}\fg_Q^2 & \dim\cO_a = 18a+12 \\
 & \dim\fr(a)=8a \\
\hspace{5mm} [\tilde A_2,2A_2,2A_2,2A_2]\hspace{3cm} 
 &\fh(a) = \fg_2, \; \fg_2,\; \fsl_2\times\fg_2, \; 2\fg_2 \hspace*{2cm}
\end{array}$$
\noindent For this case $\fh(a)=\der\AA\times\der\OO$, a product of derivation algebras. 
The grading is the doubling of the grading for $\fg_Q$.
$$\begin{array}{ll}
\fg_2\fg_Q & \dim\cO_a = 18a+18 \\
 &  \dim\fr(a)=8a+5 \\
\, [\tilde{A_2}+A_1,2A_2+A_1,2A_2+A_1,2A_2+A_1]\hspace{5mm} 
 & \fh(a) = \fsl_2, \; \fsl_2,\; 2\fsl_2, \; \fsl_2\times\fg_2 
\end{array}$$
\noindent In this case $\fh(a)=\fsl_2\times\der \AA$. Moreover, $\fg(a,0)=\fsl_2\times\CC^2
\times\fk(a)$, with the notations of the series $\fg_3$, and $\fg(a,1)=U\otimes W(a)\op V(a)$
has dimension $4a+4$, $\fg(a,2)=U\otimes V(a)\op\CC$ has dimension $4a+1$, $\fg(a,3)=U\op V(a)$
has dimension $2a+2$, $\fg(a,4)=W(a)$ and $\fg(a,5)=U$, the natural representation of $\fsl_2$. 
$$\begin{array}{ll}
\fg\fg_3 & \dim\cO_a = 18a+20 \\
 & \dim\fr(a)=7a+4 \\
 \, [C_3(a_1),A_3+A_1,\tilde{A_3}+\tilde{A_1},A_3+A_1]\qquad 
 & \fh(a) = \fsl_2, \; \fgl_2,\; 3\fsl_2, \; \fsl_2\times\fso_7
\end{array}$$
\noindent This case is similar to the previous one, since
$\fh(a)=\fsl_2\times\fl(\AA)$ and 
$\fg(a,0)=\fsl_2\times\CC\times\fk(\AA)$. But the 
induced grading is different: $\fg(a,1)=U\op U\ot V(a)$, $\fg(a,2)=V(a)\op W(a)$,
$\fg(a,3)=U\ot W(a)$, $\fg(a,4)=V(a)$,   $\fg(a,5)=U$ and $\fg(a,6)=\CC$.   

$$\begin{array}{ll}
\hspace{1cm}\fg_2^2 &  \dim\cO_a = 18a+22 \\
 & \dim\fr(a)=6a+6 \\
\hspace{1cm} [(53),F_4(a_3),D_4(a_1),D_4(a_1),D_4(a_1)]\qquad
 & \fh(a) = 0, \; 0,\; 2\CC, \; 3\fsl_2, \; \fso_8\hspace{15mm}  
\end{array}$$
\noindent Note that $\fh(a)=\ft(\AA)$, the triality algebra. The induced grading 
is the same as for the series $\fg_2$ only with indices doubled. 
$$\begin{array}{ll}
\hspace{1cm}\fg^2\fg_2^2
 & \dim\cO_a = 18a+24 \\  & \dim\fr(a)=3a+4 \\
\hspace{1cm} [(71),B_3,D_4,D_4,D_4] \hspace{3.5cm} 
 &\fh(a) = 0, \; \fsl_2,\; \fsl_3, \; \fsp_6, \; \ff_4 \hspace{2cm} 
\end{array}$$
\noindent This is the line $a=1$ of the Freudenthal square, that is $\fh(a)=\fg(\AA,\RR)$.
$$\begin{array}{ll} 
\fg\fg_3\fg_Q^2 & \dim\cO_a = 22a+20 \\
 &  \dim\fr(a)=4a+3 \\
 \, [C_3,A_5,\tilde{A_5},A_5]\hspace{3cm}
 & \fh(a) = \fsl_2, \; \fsl_2,\; 2\fsl_2, \; \fsl_2\times\fg_2\hspace{1.5cm} 
\end{array}$$
\noindent Here $\fh(a)=\fsl_2\times\der \AA$. Moreover, $\fg(a,0)=\fsl_2\times\CC^3
\times\der \AA$, and the induced grading has ten non-zero terms in positive degrees. 
$$\begin{array}{ll}
\fg_2^2 \fg_Q^2 & \dim\cO_a = 22a+22 \\
 & \dim\fr(a)=4a+4 \\
\, [F_4(a_2),E_6(a_3),E_6(a_3),E_6(a_3)] \qquad & \fh(a) = 0, \; 0,\; \fsl_2, \; \fg_2 \\
 & \\
\fg^2\fg_2^2 \fg_Q^2 & \dim\cO_a = 22a+24 \\
  & \dim\fr(a)=3a+3 \\
 \, [F_4(a_1),D_5,D_5,D_5] & \fh(a) = 0, \; \CC,\; 2\fsl_2, \; \fso_7  \\
  & \\
\fg^2\fg_2^2\fg_3^2\fg_Q^2 & \dim\cO_a = 24a+24 \\
 &  \dim\fr(a)=2a+2 \\
 \,  [F_4,E_6,E_6,E_6] & \fh(a) = 0, \; 0,\; \fsl_2, \; \fg_2 \\
\end{array}$$
\noindent We see that $\fh(a)=\der \AA$ for the two series 
$\fg_2^2 \fg_Q^2$ and $\fg^2\fg_2^2\fg_3^2\fg_Q^2$, and that 
$\fh(a)$ is given by the intermediate series $\fl(\AA)$ in the case of 
$\fg^2\fg_2^2 \fg_Q^2$. 

\medskip

\subsection{Desingularizations of orbit closures}
Given a $\fsl_2$-triple $(X,H,Y)$ in a simple complex Lie algebra $\fg$,  a 
resolution of singularities for the orbit closure $\overline{GX}$ of the
adjoint group $G$ can be obtained as follows (see \cite{pan0}): let 
$\fm=\op_{i\ge 2}\fg(i)$, 
$\fp=\op_{i\ge 0}\fg(i)$, let $P\subset G$ be the parabolic subgroup with Lie 
algebra
$\fp$. Then $\fm$ is a $P$-module, and the ``collapsing''
$$\begin{array}{ccc}
G\times_{P}\fm & \lra & \overline{GX} \subset\fg \\
\downarrow & & \\
G/P & & 
\end{array}$$
is a resolution of singularities. Here, as usual, $G\times_{P}\fm$ denotes 
the homogeneous vector bundle over the projective variety $G/P$, whose fiber at the base
point $P/P$ is the $P$-module 
$\fm$. This manifold can also be defined as the quotient of the product $G\times\fm$
by the equivalence relation $(g,m)\simeq (gp^{-1},p.m)$, where $p\in P$, so that the
product map $(g,m)\mapsto g.m\in\fg$ descends to $G\times_P\fm$. 

\smallskip Now, if $(X,H,Y)$ defines a series $\cO_a$ of nilpotent orbits in $\fg(\AA,\OO)$,
we observed that each eigenspace $\fg(a,i)$ of $ad(H)$
for the eigenvalue $i\ne 0$, so a fortiori $\fm(a)=\oplus_{i\ge 2}\fg(a,i)$, has a 
dimension which is linear in $a$. 
 This implies that the closure of $\cO_a$ is birational to a homogeneous vector 
bundle {\sl whose fiber and base are both of dimension linear in $a$}. 

Note that in most cases the orbit $\cO_a$ is even, meaning that the associated weighted
Dynkin diagram has only even weights. Such an orbit is a Richardson orbit, and the 
desingularization above of its closure is simply given by the cotangent bundle $T^*G/P$. 

\smallskip For another nice situation, consider an orbit $\overline{GX}$ corresponding 
to an $\fsl_2$ triple $(X,H,Y)$ such that $H=H_{\b}$ for some simple root $\b$. Suppose 
that $H$ defines a $5$-step grading of $\fg$, which means that the coefficient of the 
highest root $\tilde{\a}$ over $\b$ equals two. Let $P_{\b}$ denote the standard maximal 
parabolic subgroup of $G$ defined by $\b$. Consider $\tilde{\a}$ as a weight of $P_{\b}$ 
and denote by $E_{\b}(\tilde{\a})$ the associated irreducible vector bundle on $G/P_{\b}$. 
Then the desingularisation of $\overline{GX}$ is
$$\begin{array}{ccc}
E_{\b}(\tilde{\a}) & \lra & \overline{GX} \subset\fg \\
\downarrow & & \\
G/P_{\b} & & 
\end{array}$$
Recall from \cite{tits} that the adjoint variety $X_{ad}\subset\PP\fg$ is uniruled 
by the {\sl shadows} of $G/P_{\b}$, a family of homogeneous varieties parametrized
by $G/P_{\b}$. These shadows are determined pictorially by deleting $\b$ from the 
Dynkin diagram of $\fg$ with the adjoint marking (when the adjoint representation 
is fundamental, this just means that we mark the node of the corresponding fundamental
weight). Then the projectivization of $\overline{GX}$ is the union of the linear
spans of the shadows, and the vector bundle $E_{\b}(\tilde{\a})$ is the family of 
the associated vector subspaces of $\fg$. (Special cases of this were observed in 
\cite{LMmagic}.) 
This phenomenon occurs uniformly for the series $\fg_Q$.


\subsection{Rational points}
A nilpotent orbit $\cO\simeq G/K\subset\fg$ is defined over $\FF_q$ for $q$ large enough, 
and the number of its $\FF_q$-points is a polynomial function of $q$ ( 
\cite{bp}, Theorem 1.a). 
We can deduce this polynomial function from the data gathered in \cite{carter}.
Indeed, if $K$ is connected this number is equal to $|G(\FF_q)|/|K(\FF_q)|$ (see 
\cite{bp}, Theorem 1.c), and can be deduced from the formulas in \cite{carter}, pp. 75-76, 
and the data for $K$ gathered above. When the group $K$ is not connected, which may happen 
in some cases, the formulas below hold for the quotients $|G(\FF_q)|/|K(\FF_q)|$.

For each of our series $\cO_a$ of nilpotent orbits, we express the resulting  polynomial 
as a rational 
function involving only terms of the form $q^{\ell}-1$, where $\ell$ is some linear function 
of $a$, from a very limited list. 

We begin with the biggest series of orbits, whose label is $\fg^2\fg_2^2\fg_3^2\fg_Q^2$.
The number of $\FF_q$-points on these orbits is 
{\small $$Z_{\fg^2\fg_2^2\fg_3^2\fg_Q^2}(q)=q^{11a+8}\frac{(q^a-1)(q^{3a/2}-1)(q^{3a/2+2}-1)
(q^{2a+2}-1)(q^{2a+4}-1)(q^{5a/2+4}-1)(q^{3a+6}-1)}{(q^{a/2+2}-1)}.$$ }
For the other series, the corresponding functions are simple quotients 
$Z_{\cO}(q)=Z_{\fg^2\fg_2^2\fg_3^2\fg_Q^2}(q)/Y_{\cO}(q)$, with denominators
given by the following table:

{\small 
\begin{eqnarray}
\nonumber 
Y_{\fg}(q) & = & q^{11a+8}(q^{a}-1)(q^{a+2}-1)(q^{3a/2}-1)(q^{3a/2+2}-1)(q^{2a+2}-1),  \\
\nonumber 
Y_{\fg_Q}(q) & = & q^{21a/2+6}(q^{a/2}-1)(q^{a}-1)(q^{a+2}-1)(q^{a+4}-1),  \\
\nonumber 
Y_{\fg_2}(q) & = & q^{19a/2+6}(q^{2}-1)(q^{a}-1)(q^{3a/2}-1),  \\
\nonumber 
Y_{\fg^2}(q) & = & q^{8a+4}(q^{a/2+1}-1)(q^{a}-1)(q^{a+1}-1)(q^{3a/2}-1),  \\
\nonumber 
Y_{\fg_3}(q) & = & q^{15a/2+4}(q^{2}-1)(q^{a/2}-1),  \\
\nonumber 
Y_{\fg^2\fg_Q}(q) & = & q^{11a/2+2}(q^{a/2}-1)(q^{a}-1)(q^{a+2}-1),  \\
\nonumber 
Y_{\fg_Q^2}(q) & = & q^{6a+4}(q^{2}-1)(q^{6}-1),  \\
\nonumber 
Y_{\fg_2\fg_Q}(q) & = & q^{6a+4}(q^{2}-1),  \\
\nonumber 
Y_{\fg\fg_3}(q) & = & q^{11a/2+2}(q^{2}-1)(q^{a/2}-1),  \\
\nonumber 
Y_{\fg_2^2}(q) & = & q^{5a+2}(q^{a/2}-1)^2,  \\
\nonumber 
Y_{\fg^2\fg_2^2}(q) & = & q^{7a/2}(q^{a}-1)(q^{3a/2}-1),  \\
\nonumber 
Y_{\fg\fg_3\fg_Q^2}(q) & = & q^{2a+2}(q^{2}-1),  \\
\nonumber 
Y_{\fg_2^2\fg_Q^2}(q) & = & q^{2a+2},  \\
\nonumber 
Y_{\fg^2\fg_2^2\fg_Q^2}(q) & = & q^{3a/2}(q^{a/2}-1).
\end{eqnarray}  }

In particular, the number of $\FF_q$-points on the series of minimal nilpotent
orbits is 
{\small $$ Z_{\fg}(q)=\frac{(q^{2a+4}-1)(q^{5a/2+4}-1)(q^{3a+6}-1)}
{(q^{a/2+2}-1)(q^{a+2}-1)}.$$ }

\bigskip
\subsection{Unipotent characters}

The Springer correspondence uses local systems on nilpotent orbits 
to define representations of Weyl groups, which themselves are in 
natural correspondence with unipotent characters of finite groups of 
the corresponding Lie type. 
In this section we show that the unipotent characters corresponding to our
series of nilpotent orbits in the exceptional Lie algebras  
are accordingly organized into series.  This can be seen on the
polynomials giving the degrees of these characters, once we write
these polynomials as rational functions. More precisely, we are able to write 
these functions as products of factors of type $q^e-1$, or inverses of such factors, 
with $e$ a linear function of $a$. This   striking fact only holds 
for $a=2, 4, 8$. A theoretical explanation would be most welcome. Also 
it would be interesting to understand what really happens when $a=1$, that
is, when $\fe_6$ is folded into $\ff_4$. 

\smallskip
Note that the fundamental groups of the nilpotent orbits in our series 
are well-behaved: they are constant in each series, either trivial
or equal to $\ZZ_2$, in which case we get two series of unipotent 
characters. Actually there is one exception to this: in the series
labeled $\fg_Q^2$, the nilpotent orbits of $\fe_6$ and $\fe_7$ are 
simply connected, but that of $\fe_8$ has fundamental group $\ZZ_2$. 

\smallskip
The following data are again transcriptions of the formulas gathered in \cite{carter}
(pp. 480-488) for the degrees of unipotent characters. Note that in this reference, these degrees
are given as products of cyclotomic polynomials, a form in which the regularities
that we observed are far from visible. Some work is needed to put these formulas 
into the form that follows. 
Note that only a small family of linear functions are involved in these
formulas. 
Note also that many simplifications may occur in each degree, but in different ways. 

We let $N$ denote the number of positive roots.  

\bigskip\noindent $\fg$ : The degree of the associated unipotent character is

{\small \[ q^{N-3a-5}\frac{(q^{2a+4}-1)(q^{5a/2+4}-1)}
{(q^{a/2+2}-1)(q^{a+2}-1)} \] }

\bigskip\noindent $\fg_Q$ : The degree of the associated unipotent character is
{\small \[ q^{N-5a-6}
\frac{(q^{3a/2}-1)(q^{3a/2+2}-1)(q^{2a+4}-1)(q^{3a+6}-1)}
{(q^{a/2}-1)(q^{a/2+2}-1)(q^{a+2}-1)(q^{a+4}-1)} \] }

\bigskip
\noindent $\fg_2$ : The degree of the associated unipotent character is

{\small \[ \frac{1}{2} q^{N-6a-9}
\frac{(q^{a/2+1}-1)(q^{a+1}-1)(q^{3a/2+2}-1)(q^{2a+4}-1)(q^{5a/2+4}-1)(q^{3a+6}-1)}
{(q-1)(q^{a/2+2}-1)^2(q^{a/2+1}-1)(q^{a+2}-1)^2(q^{3a/2+3}-1)} \] }

\bigskip\noindent $\fg^2$ : The degrees of the two associated unipotent characters are

{\small \[ \frac{1}{2} q^{N-6a-9}
\frac{(q^{a+2}-1)(q^{3a/2+2}-1)(q^{2a+2}-1)(q^{5a/2+4}-1)(q^{3a+6}-1)}
{(q-1)(q^{a/2+1}-1)(q^{a+1}-1)(q^{a+4}-1)(q^{3a/2+3}-1)} \]
\[ \frac{1}{2} q^{N-6a-9}
\frac{(q-1)(q^{3a/2+2}-1)(q^{3a/2+3}-1)(q^{2a+2}-1)(q^{2a+4}-1)(q^{5a/2+4}-1)}
{(q^2-1)(q^{a/2+1}-1)(q^{a/2+2}-1)^2(q^{a+1}-1)(q^{a+2}-1)} \] }

\bigskip\noindent $\fg^2\fg_Q$ : The degree of the associated unipotent character is

{\small \[ q^{N-8a-10}
\frac{(q^{3a/2}-1)(q^{3a/2+2}-1)(q^{2a+2}-1)(q^{5a/2+4}-1)(q^{3a+6}-1)}
{(q^2-1)(q^{a/2}-1)(q^{a/2+2}-1)(q^{a/2+4}-1)(q^{a+2}-1)} \] }

\bigskip \noindent $\fg_2^2$ :  The degrees of the two associated 
unipotent characters are $q^{N-9a-11}$ times 

{\small \[ 
\frac{(q^{a/2+1}-1)^3(q^{a}-1)(q^{3a/2}-1)(q^{3a/2+2}-1)
(q^{2a+2}-1)(q^{2a+4}-1)(q^{5a/2+4}-1)(q^{3a+6}-1)}
{6(q-1)^2(q^2-q+1)(q^{a/2}-1)^2(q^{a/2+2}-1)^3(q^{a+2}-1)^2(q^{3a/2+3}-1)} \] 
\[ 
\frac{(q^a-1)(q^{3a/2}-1)(q^{3a/2+2}-1)(q^{2a+2}-1)(q^{2a+4}-1)
(q^{5a/2+4}-1)(q^{3a+6}-1)}
{3(q^2-1)^2(q^{a/2}-1)^2(q^{a+2}-1)^2(q^{3a/2+6}-1)} \] }

\pagebreak
\noindent $\fg^2\fg_2^2$ :
The degree of the associated unipotent character is

{\small \[ q^{N-9a-12}
\frac{(q^{3a/2+2}-1)(q^{2a+2}-1)(q^{2a+4}-1)(q^{5a/2+4}-1)(q^{3a+6}-1)}
{(q^2-1)(q^6-1)(q^{a/2+2}-1)(q^{a/2+4}-1)(q^{a+4}-1)} \] }

\bigskip
\noindent $\fg_Q^2$ :
Here there is a problem: there are two associated characters for $E_8$, 
but only one for $E_6$ and $E_7$. Nevertheless, let 

{\small  \[ \phi_a(q) = q^{N-9a-6}
\frac{(q^a-1)(q^{3a/2+2}-1)(q^{2a+4}-1)(q^{5a/2+4}-1)(q^{3a+6}-1)}
{(q^2-1)^2(q^6-1)(q^{a/2+2}-1)(q^{a/2+4}-1)} \] }

The degrees of the unipotent characters attached to this series for $E_6$ and $E_7$
are $\phi_2(q)$ and $\phi_4(q)$, while the two characters for $E_8$ 
have their degrees given by 

{\small $$\phi_{8,\epsilon}(q)=\frac{1}{2}\frac{q^9-\epsilon}{q^3-\epsilon}
\frac{q-\epsilon}{q^7-\epsilon}\phi_8(q), \qquad \epsilon=\pm 1.$$ }

\bigskip\noindent $\fg_3$ : The degree of the associated unipotent character is
{\small \[ q^{N-8a-9}
\frac{(q^{a/2+4}-1)(q^{2a-2}-1)(q^{2a+4}-1)
(q^{5a/2+4}-1)(q^{3a+6}-1)}
{(q^2-1)(q^6-1)(q^{a/2}-1)(q^{a/2+2}-1)(q^{a+2}-1)} \] }

\bigskip\noindent $\fg_2\fg_Q$ :
The degree of the associated unipotent character is
 
{\small \[ \frac{1}{3}q^{N-9a-11}
\frac{(q^{a/2}-1)(q^a-1)(q^{3a/2+2}-1)(q^{2a+2}-1)(q^{2a+4}-1)
(q^{5a/2+4}-1)(q^{3a+6}-1)}
{(q^2-1)^2(q^{a/2+2}-1)^3(q^{a+2}-1)^2} \] }

\bigskip\noindent $\fg\fg_3$ :
The degree of the associated unipotent character is

{\small \[ \frac{1}{2}q^{N-9a-11}
\frac{(q^{3a/2}-1)(q^{3a/2+2}-1)(q^{2a+2}-1)(q^{2a+4}-1)
(q^{5a/2+4}-1)(q^{3a+6}-1)}
{(q-1)(q^3-1)(q^{a/2+1}-1)(q^{a/2+2}-1)(q^{a+4}-1)(q^{3a/2+3}-1)} \] }

\bigskip\noindent $\fg\fg_3\fg_Q^2$ :
The degree of the associated unipotent character is
$\frac{1}{2}q^{N-11a-11}$ times 

{\small \[ 
\frac{(q^a-1)(q^{3a/2}-1)(q^{3a/2+2}-1)(q^{2a+2}-1)(q^{2a+4}-1)
(q^{5a/2+4}-1)(q^{3a+6}-1)}
{(q-1)(q^3-1)(q^4-1)(q^{a/2+2}-1)(q^{a/2+3}-1)(q^{a/2+5}-1)(q^{a+4}-1)} \] }

\bigskip\noindent $\fg_2^2\fg_Q^2$ : The degree of the associated character is 
$\frac{1}{2}q^{N-11a-11}$ times 

{\small \[ 
\frac{(q^{a/2+3}-1)(q^a-1)(q^{3a/2}-1)(q^{3a/2+2}-1)(q^{2a+2}-1)
(q^{2a+4}-1)(q^{5a/2+4}-1)(q^{3a+6}-1)}
{(q-1)(q^2-1)^2(q^3+1)(q^{a/2+2}-1)^3(q^{a/2+5}-1)(q^{a+6}-1)} \] }

\bigskip\noindent $\fg^2\fg_2^2\fg_Q^2$ : 
The degree of the associated unipotent character is 
$q^{N-11a-12}$ times 

{\small \[ 
\frac{(q^{a/2+4}-1)(q^{2a-2}-1)(q^{2a+2}-1)
(q^{2a+4}-1)(q^{5a/2+4}-1)(q^{3a+6}-1)}
{(q^2-1)(q^4-1)(q^6-1)^2(q^{a/2}-1)(q^{a/2+8}-1)} \]
}

\bigskip\noindent $\fg^2\fg_2^2\fg_3^2\fg_Q^2$ : 
The degree of the associated unipotent character is 
$q^{N-12a-12}$ times 

{\small \[ 
\frac{(q^a-1)(q^{3a/2}-1)(q^{3a/2+2}-1)(q^{2a+2}-1)
(q^{2a+4}-1)(q^{5a/2+4}-1)(q^{3a+6}-1)}
{(q^2-1)(q^6-1)(q^8-1)(q^{12}-1)(q^{a/2+2}-1)(q^{a/2+4}-1)(q^{a/2+8}-1)} \]
}


\subsection{Series of type $E_6$}

We now examine how the five remaining nilpotent orbits in $\fe_6$ propagate to
orbits in $\fe_7$ and $\fe_8$. They are associated
to $\fsl_2$-triples $(X,H,Y)$ for which the semi-simple 
element $H$ can be chosen to belong to $\ft(\OO)=\fso_8$, hence can again be encoded
by a label $\fg^p\fg_2^q\fg_3^r\fg_Q^s$.

\smallskip
The degrees of the associated unipotent characters do not behave as well as in 
the series coming from $\ff_4$. A first difficulty is that in  each case, 
there are two associated characters in type $E_7$ and $E_8$, but only one 
in type $E_6$. We already encountered a similar phenomenon for the series 
$\fg_Q^2=[A_2,2A_2,2A_2,2A_2]$, where the degrees of the two 
unipotent characters were closely related. This is again true for the 
series of type $E_6$, and an a priori explanation  would be welcome.

\smallskip
For each series of orbits, we provide the label used in \cite{carter}, the dimension
of the orbits $\cO_a$ and of the unipotent radical $\fr(a)$ of the generic centralizers, which
again are both linear functions in $a$, and the reductive parts $\fh(a)$ of these centralizers.

\bigskip\noindent $\fg\fg_Q,\quad $ $A_2+A_1$, $\qquad\dim\cO_a = 15a+16$, 
$\qquad\dim\fr(a)=9a+5, \quad \fh(a) = \fgl_3, \; \fgl_4,\; \fsl_6.$
\medskip

\noindent The degree of the associated unipotent character in type $E_6$ is 
{\small $$\deg\phi_{64,13}=q^{13}\frac{(q^6-1)(q^8-1)(q^{12}-1)}{(q-1)(q^3-1)(q^3-1)}.$$}
In type $E_7$ the degrees of the two unipotent characters are
{\small $$\deg\phi_{120,25}=\frac{1}{2}q^{25}\frac{(q^8-1)(q^{10}-1)(q^{12}-1)(
q^{18}-1)}{(q-1)(q^3-1)(q^4-1)(q^6-1)}
\times\frac{(q^3+1)(q^7+1)}{(q^4+1)(q^6+1)},$$
$$\deg\phi_{105,26}=\frac{1}{2}q^{25}\frac{(q^8-1)(q^{10}-1)(q^{12}-1)(
q^{18}-1)}{(q-1)(q^3-1)(q^4-1)(q^6-1)}
\times\frac{(q^3-1)(q^7-1)}{(q^4-1)(q^6-1)}.$$}
 
\medskip\noindent
In type $E_8$ the degrees of the two unipotent characters are
{\small $$\deg\phi_{210,52}=\frac{1}{2}q^{52}\frac{(q^{14}-1)(q^{18}-1)(q^{20}-1)(
q^{30}-1)}{(q^3-1)(q^4-1)(q^5-1)(q^6-1)}
\times\frac{(q^4+1)(q^{12}+1)}{(q^7+1)(q^9+1)},$$
$$\deg\phi_{160,55}=\frac{1}{2}q^{52}\frac{(q^{14}-1)(q^{18}-1)(q^{20}-1)(
q^{30}-1)}{(q^3-1)(q^4-1)(q^5-1)(q^6-1)}
\times\frac{(q^4-1)(q^{12}-1)}{(q^7-1)(q^9-1)}.$$}
\medskip

Let us introduce the following rational function, which is close to those
we already met, except for the appearance of an $a/4$ in the exponents of $q$:

{\small 
\[ \psi_{\fg\fg_Q}(q)=\frac{1}{2}
\frac{(q^2-1)(q^{5a/4-2}-1)(q^{3a/2+2}-1)(q^{2a+2}-1)(q^{2a+4}-1)(q^{3a+6}-1)}
{(q^3-1)(q^{a/4-1}-1)(q^{a/2}-1)(q^{a/2+1}-1)(q^{a/2+2}-1)(q^{a/2+4}-1)}. \] }
Then we can write the degrees above as
{\small $$ q^{N-15a/2-8}\psi_{\fg\fg_Q}(q)
\frac{(q^{a/4+2}-1)(q^{5a/4+2}-1)}{(q^{3a/4+1}-1)(q^{3a/4+3}-1)}
\quad {\rm and} \quad
q^{N-15a/2-8}\psi_{\fg\fg_Q}(q)
\frac{(q^{a/4+2}+1)(q^{5a/4+2}+1)}{(q^{3a/4+1}+1)(q^{3a/4+3}+1)}.$$}

These formulas have several intriguing features. They are obviously 
closely related one to the other. 
For $a=1$, the non integer 
exponents cancel out. Moreover, the second part of this expression gives $1$ 
for $a=1$, hence the same rational expression with coefficient one half : in 
fact, there is only one character in this case, whose degree is given by the
sum of these two equal contributions. What kind of group theoretic explanation
could this phenomenon have ? 

\smallskip\noindent In  this series, the number of $\FF_q$-points is given by 
{\small 
$$q^{3a+4}\frac{(q^{5a/4-2}-1)(q^{3a/2}-1)(q^{3a/2+2}-1)
(q^{2a+2}-1)(q^{2a+4}-1)(q^{5a/2+4}-1)(q^{3a+6}-1)}
{(q^3-1)(q^{a/4}-1)(q^{a/2}-1)(q^{a/2+1}-1)(q^{a/2+2}-1)}.$$}

\bigskip\noindent $\fg^2\fg_Q^2$, $\quad A_4$, $\qquad\dim\cO_a = 20a+20,
\qquad \dim\fr(a)=5a+4, \quad \fh(a) = \fgl_2,\; \fgl_3, \; \fsl_5.$
\medskip

\noindent Here we have one unipotent character $\phi_{81,6}$ in type $E_6$, two in type 
$E_7$, $\phi_{420,13}$ and $\phi_{336,14}$, and again two in type $E_8$, 
$\phi_{2268,30}$ and $\phi_{1296,33}$. Their degrees are given by the following 
expressions, with the same phenomenon for $a=1$ as in the previous case:

{\small 
\[ q^{N-10a-10} \psi_{\fg^2\fg_Q^2}(q)
\frac{(q^{2}-1)(q^{a+2}-1)}{(q^{a/2+1}-1)(q^{a/2+3}-1)} 
\quad {\rm and} \quad 
q^{N-10a-10} \psi_{\fg^2\fg_Q^2}(q)
\frac{(q^{2}+1)(q^{a+2}+1)}{(q^{a/2+1}+1)(q^{a/2+3}+1)},  \] 

\[ \psi_{\fg^2\fg_Q^2}(q)=\frac{1}{2}
\frac{(q^{3a/4-1}-1)(q^{3a/4}-1)(q^{a+4}-1)(q^{2a-2}-1)(q^{2a+2}-1)(q^{5a/2+4}-1)(q^{3a+6}-1)}
{(q^4-1)(q^6-1)(q^{a/4}-1)(q^{a/4+1}-1)(q^{a/2}-1)(q^{a/2+1}-1)^2}. \] }

\medskip\noindent In  this series, the number of $\FF_q$-points is given by 
{\small 
$$q^{15a/2+6}\frac{(q^2-1)(q^{5a/4-2}-1)(q^{3a/2}-1)(q^{3a/2+2}-1)(q^{2a+2}-1)(q^{2a+4}-1)
(q^{5a/2+4}-1)(q^{3a+6}-1)}
{(q^3-1)(q^{a/4}-1)(q^{a/2}-1)(q^{a/2+1}-1)}.$$}
\medskip

\bigskip\noindent $\fg\fg_3\fg_Q$, $\quad A_4+A_1$, $\qquad\dim\cO_a = 21a+20, 
\qquad \dim\fr(a)=6a+3, \quad \fh(a) =  \CC\; \CC^2, \; \fgl_3.$
\medskip

\noindent Here we have one unipotent character $\phi_{60,5}$ in type $E_6$, two in type 
$E_7$, $\phi_{512,11}$ and $\phi_{512,12}$, and again two in type $E_8$, 
$\phi_{4096,26}$ and $\phi_{4096,27}$. Their degrees are given by 
$$q^{N-21a/2-10} \psi_{\fg\fg_3\fg_Q}(q),$$
(except for $a=1$ where the degree of the unique character is twice 
this quantity), with 
{\small
\[ \psi_{\fg\fg_3\fg_Q}(q)=\frac{1}{2}
\frac{(q^2-1)(q^{5a/4-2}-1)(q^{3a/2}-1)(q^{3a/2+2}-1)(q^{2a+2}-1)(q^{2a+4}-1)(q^{5a/2+4}-1)(q^{3a+6}-1)}
{(q-1)(q^3-1)^2(q^{a/4}-1)(q^{a/2+1}-1)(q^{a/2+3}-1)(q^{a/2+5}-1)(q^{3a/2+3}-1)}. \]  }

\medskip\noindent In  this series, the number of $\FF_q$-points is given by 
{\small 
$$q^{15a/2+6}\frac{(q^2-1)(q^{5a/4-2}-1)(q^{3a/2}-1)(q^{3a/2+2}-1)(q^{2a+2}-1)(q^{2a+4}-1)
(q^{5a/2+4}-1)(q^{3a+6}-1)}
{(q-1)(q^3-1)(q^{a/4}-1)}.$$}

\bigskip\noindent $\fg^2\fg_3\fg_Q$, $\quad D_5(a_1)$, $\qquad\dim\cO_a = 21a+22, 
\qquad \dim\fr(a)=5a+3, \quad \fh(a) = \CC, \; \fgl_2, \; \fsl_4.$

\medskip\noindent
Here we have one unipotent character $\phi_{64,4}$ in type $E_6$, two in type 
$E_7$, $\phi_{420,10}$ and $\phi_{336,11}$, and again two in type $E_8$, 
$\phi_{2800,25}$ and $\phi_{2100,28}$. Their degrees are given by 
$$q^{N-21a/2-11} \psi_{\fg^2\fg_3\fg_Q}(q)
\quad {\rm and} \quad 
q^{N-21a/2-11} \psi_{\fg^2\fg_3\fg_Q}'(q),$$
(except for $a=1$ where the degree of the unique character is the sum of  
these two -- equal in this case --  quantities), with 

{\small 
\[ \psi_{\fg^2\fg_3\fg_Q}(q)=\frac{1}{2}
\frac{(q^{a/4+4}-1)(q^{3a/4}-1)(q^{5a/4-2}-1)(q^{3a/2+2}-1)(q^{2a+2}-1)
(q^{2a+4}-1)(q^{5a/2+4}-1)(q^{3a+6}-1)}
{(q^3-1)^2(q^{a/4}-1)(q^{a/4+1}-1)(q^{a/2}-1)(q^{a/2+4}-1)(q^{a/2+8}-1)(q^{3a/4+3}-1)}, \]

\[ \psi_{\fg^2\fg_3\fg_Q}'(q)=\frac{1}{2}
\frac{(q^{3a/4-1}-1)(q^{5a/4-1}-1)(q^{3a/2}-1)(q^{3a/2+2}-1)(q^{2a+2}-1)
(q^{2a+4}-1)(q^{5a/2+4}-1)(q^{3a+6}-1)}
{(q^3-1)(q^5-1)(q^{a/4}-1)(q^{a/2}-1)(q^{a/2+2}-1)^2(q^{3a/4}-1)(q^{3a/2+6}-1)}. \]   }

\medskip\noindent In  this series, the number of $\FF_q$-points is given by 
{\small 
$$q^{8a+7}\frac{(q^2-1)(q^{5a/4-2}-1)(q^{3a/2}-1)(q^{3a/2+2}-1)(q^{2a+2}-1)(q^{2a+4}-1)
(q^{5a/2+4}-1)(q^{3a+6}-1)}
{(q^3-1)(q^{a/4}-1)(q^{a/2}-1)}.$$}

\pagebreak

\noindent $\fg^2\fg_3^2\fg_Q^2,\quad$ $E_6(a_1)$, 
$\qquad\dim\cO_a = 24a+22, 
\qquad \dim\fr(a)=3a+2, \quad \fh(a) = 0, \; 0,\; \fsl_3.$

\medskip\noindent
Here again we have one unipotent character $\phi_{6,1}$ in type $E_6$, two in type 
$E_7$, $\phi_{120,4}$ and $\phi_{105,5}$, and again two in type $E_8$, 
$\phi_{2800,13}$ and $\phi_{2100,16}$. Their degrees are given by 
$$q^{N-12a-11} \psi_{\fg^2\fg_3^2\fg_Q^2}(q)
\quad {\rm and} \quad 
q^{N-21a/2-11} \psi_{\fg^2\fg_3^2\fg_Q^2}'(q),$$
(except for $a=1$ where the degree of the unique character is the sum of  
these two -- equal in this case --  quantities), with 

{\small
\[ \psi_{\fg^2\fg_3^2\fg_Q^2}(q)=\frac{1}{2}
\frac{(q^{a/2+2}-1)(q^{3a/4}-1)(q^{3a/4}-1)(q^{2a-2}-1)(q^{2a+2}-1)
(q^{2a+4}-1)(q^{5a/2+4}-1)(q^{3a+6}-1)}
{(q^3-1)^2(q^4-1)(q^{12}-1)(q^{a/4}-1)(q^{a/4+1}-1)(q^{a/2+1}-1)(q^{a/2+5}-1)}, \]

\[ \psi_{\fg^2\fg_3^2\fg_Q^2}'(q)=\frac{1}{2}
\frac{(q^{a/2+5}-1)(q^{5a/4-2}-1)(q^{3a/2}-1)(q^{3a/2+2}-1)(q^{2a+2}-1)
(q^{2a+4}-1)(q^{5a/2+4}-1)(q^{3a+6}-1)}
{(q^3-1)(q^4-1)(q^6-1)^2(q^{a/4}-1)(q^{3a/2}-1)(q^{3a/2+2}-1)(q^{a+10}-1)}. \] }

\medskip\noindent In  this series, the number of $\FF_q$-points is given by 
{\small 
$$q^{21a/2+7}\frac{(q^2-1)(q^{5a/4-2}-1)(q^{3a/2}-1)(q^{3a/2+2}-1)(q^{2a+2}-1)(q^{2a+4}-1)
(q^{5a/2+4}-1)(q^{3a+6}-1)}
{(q^3-1)(q^{a/4}-1)}.$$}

\medskip
This accounts for all nilpotent orbits in $\fe_6$, about one half of those in $\fe_7$ 
and a little less than one third of those in $\fe_8$. 

\section{Series for the other rows of Freudenthal square}

The exceptional series of Lie algebras is the fourth line $\fg(\AA,\OO)$
in the magic square of Freudenthal, and we just saw how this allows us to 
organize their nilpotent orbits into series. 

In this section we briefly discuss the other three lines of Freudenthal
square and their nilpotent orbits. 

\subsection{The subexceptional series $\fg(\AA,\HH)$}
Here the Lie algebras $\fg$, and the number of positive roots $N$, 
parametrized by $a$ are: 
\[ \begin{array}{ccccc}
a & 1 & 2 & 4 & 8 \\
\fg  & \fsp_6 & \fsl_6 & \fso_{12} & \fe_7 \\
N & 9 & 15 & 30 & 63
\end{array} \]

The nilpotent orbits of $\fso_{12}$ are parametrized by pairs of 
partitions $(\alpha,\beta)$ such that $2|\alpha|+|\beta|=12$
and $\beta$ has distinct parts. The nilpotent orbits of $\fsl_6$
are parametrized by partitions of six. The nilpotent orbits of
$\fsp_6$ are parametrized by pairs of partitions $(\alpha,\beta)$
with $|\alpha|+|\beta|=3$ where $\beta$ has distinct parts (see \cite{carter}).

Given a nilpotent orbit $(\alpha,\beta)$ of $\fsp_6$ the elementary divisors
are given by repeating each part of $\alpha$ twice and doubling each part
of $\beta$. By ordering these elementary divisors we get a partition
$\lambda$ with $|\lambda|=6$ which corresponds to a nilpotent orbit
of $\fsl_6$. Given a nilpotent orbit $\lambda$ of $\fsl_6$ we can take the
pair of partitions $(\lambda,\emptyset)$ which corresponds to a nilpotent
orbit of $\fso_{12}$. These constructions give the first three terms of each
series below. 

For the subexceptional series we have three preferred representations, 
$\fg, \fg_Q=V_2, \fg_{\BA\BP^2}=V$ in the notations of \cite{LMtrial}. 
The highest weights of these representations are as follows:

$$\begin{array}{rcccc}
 & \fsp_6 & \fsl_6\;\; & \fso_{12} & \fe_7 \\
\o(\fg) & 
\setlength{\unitlength}{2mm}
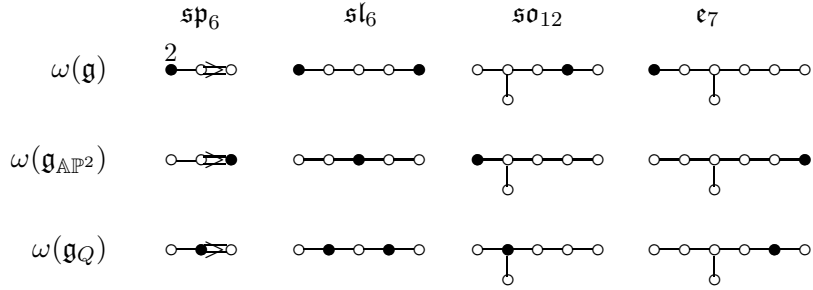
\begin{picture}(7,3)(-2,0)
\multiput(0,0)(2,0){3}{$\circ$}
\put(0,0){$\bullet$}
\put(0.5,.5){\line(1,0){1.6}}
\put(0,1.1){{\small $2$}} 
\multiput(2.6,.2)(0,.5){2}{\line(1,0){1.5}} 
\put(2.6,0){$>$}
\end{picture}\quad & 
\setlength{\unitlength}{2mm}
\begin{picture}(8,3)(0.2,0)
\multiput(0,0)(2,0){5}{$\circ$}
\multiput(0.8,.5)(2,0){4}{\line(1,0){1.3}} 
\put(0,0){$\bullet$}\put(8,0){$\bullet$}
\end{picture}\quad & 
\setlength{\unitlength}{2mm}
\begin{picture}(9,3)(0,0)
\multiput(0,0)(2,0){5}{$\circ$}
\put(6,0){$\bullet$}
\multiput(0.8,.5)(2,0){4}{\line(1,0){1.3}} 
\put(2,-2){$\circ$}
\put(2.45,-1.3){\line(0,1){1.5}}
\end{picture} & 
\setlength{\unitlength}{2mm}
\begin{picture}(10,3)(-1,0)
\multiput(0,0)(2,0){6}{$\circ$}
\put(0,0){$\bullet$}
\multiput(0.8,.5)(2,0){5}{\line(1,0){1.3}} 
\put(4,-2){$\circ$}
\put(4.45,-1.3){\line(0,1){1.4}}
\end{picture} \\
 & & & & \\
\o(\fg_{\AA\PP^2}) &
\setlength{\unitlength}{2mm}
\begin{picture}(7,3)(-2,0)
\multiput(0,0)(2,0){3}{$\circ$}
\put(4,0){$\bullet$}
\put(0.8,.4){\line(1,0){1.3}}
\multiput(2.6,.2)(0,.5){2}{\line(1,0){1.5}} 
\put(2.6,0){$>$}
\end{picture}\quad & 
\setlength{\unitlength}{2mm}
\begin{picture}(8,3)(0.2,0)
\multiput(0,0)(2,0){5}{$\circ$}
\multiput(0.8,.5)(2,0){4}{\line(1,0){1.3}} 
\put(4,0){$\bullet$}
\end{picture}\quad & 
\setlength{\unitlength}{2mm}
\begin{picture}(9,3)(0,0)
\multiput(0,0)(2,0){5}{$\circ$}
\put(0,0){$\bullet$}
\multiput(0.8,.5)(2,0){4}{\line(1,0){1.3}} 
\put(2,-2){$\circ$}
\put(2.45,-1.3){\line(0,1){1.3}}
\end{picture} & 
\setlength{\unitlength}{2mm}
\begin{picture}(10,3)(-1,0)
\multiput(0,0)(2,0){6}{$\circ$}
\put(10,0){$\bullet$}
\multiput(0.8,.5)(2,0){5}{\line(1,0){1.3}} 
\put(4,-2){$\circ$}
\put(4.45,-1.3){\line(0,1){1.3}}
\end{picture} \\
 & & & & \\
\o(\fg_Q) & 
\setlength{\unitlength}{2mm}
\begin{picture}(7,3)(-2,0)
\multiput(0,0)(2,0){3}{$\circ$}
\put(2,0){$\bullet$}
\put(0.8,.5){\line(1,0){1.3}}
\multiput(2.7,.3)(0,.5){2}{\line(1,0){1.4}} 
\put(2.6,0){$>$}
\end{picture}\quad & 
\setlength{\unitlength}{2mm}
\begin{picture}(8,3)(0.2,0)
\multiput(0,0)(2,0){5}{$\circ$}
\multiput(0.8,.5)(2,0){4}{\line(1,0){1.3}} 
\put(2,0){$\bullet$}\put(6,0){$\bullet$}
\end{picture} \quad & 
\setlength{\unitlength}{2mm}
\begin{picture}(9,3)(0,0)
\multiput(0,0)(2,0){5}{$\circ$}
\put(2,0){$\bullet$}
\multiput(0.8,.5)(2,0){4}{\line(1,0){1.3}} 
\put(2,-2){$\circ$}
\put(2.45,-1.3){\line(0,1){1.3}}
\end{picture} & 
\setlength{\unitlength}{2mm}
\begin{picture}(10,3)(-1,0)
\multiput(0,0)(2,0){6}{$\circ$}
\put(8,0){$\bullet$}
\multiput(0.8,.5)(2,0){5}{\line(1,0){1.3}} 
\put(4,-2){$\circ$}
\put(4.45,-1.3){\line(0,1){1.3}}
\end{picture} 

\end{array}$$

\bigskip
We obtain the following series:

$$\begin{array}{lll}
\fg\qquad & [(11,1),(21^4),(21^4,-),A_1]\qquad & \dim\cO_a = 4a+2, \\
 & & \dim\fr (a)  = 4a+1, \\
 & & \fh (a)= \fso_5, \fsl_4(\fso_6), \fsl_2\times\fso_8, \fso_{12}, \\
 & & \\
\fg_Q & [(21,-),(2211),(2211,-),2A_1] & \dim\cO_a =6a+4, \\
 & & \dim  \fr (a) =5a+2, \\
 & & \fh(a)=\fgl_2 , 2\fsl_2 , \fsl_2\times\fso_5 , \fsl_2\times\fso_9=\fsl_2\times\fso_{a+1}, \\
 & & \\
\fg_{\AA\PP^2} & [(2,1),(222),(222,-),3A_1] & \dim\cO_a =6a+6, \\
 & & \dim  \fr (a)  =3a+3, \\
 & & \fh (a)= \fsl_2 , \fsl_3 , \fsp_6 , \ff_4 =\fg(\AA,\RR), \\
 & & \\
\fg_Q^2 & [(3,-),(33),(33,-),2A_2] & \dim\cO_a =10a+4, \\
 & & \dim  \fr (a) =4a, \\
 & & \fh(a) = \fsl_2 , \fsl_2 , 2\fsl_2 , \fsl_2\times\fg_2 = \fsl_2\times\der \AA, \\
 & & \\
\fg_{\AA\PP^2}^2\fg_Q & [(1,2),(411),(411,-),A_3] & \dim\cO_a =10a+4, \\
 & & \dim  \fr (a) =3a+1, \\
 & & \fh (a)=\fsl_2 , \fgl_2 , 3\fsl_2 , \fsl_2\times\fso_7=
\fsl_2\times\fl(\AA), \\
 & & \\
\fg_{\AA\PP^2}\fg & [(-,21),(42),(42,-),(A_3+A_1)^{\prime\prime}] & \dim\cO_a =10a+6, \\
 & & \dim  \fr (a) =3a+2, \\
 & & \fh (a)=0, \CC , 2\fsl_2 , \fso_7=\fl(\AA), \\
 & & \\
\fg^2\fg_{\AA\PP^2}^2\fg_Q^2 & [(-,3),(6),(6,-),A_5] & \dim\cO_a =12a+6, \\
 & & \dim  \fr (a) =2a+1, \\ 
 & & \fh (a)= 0 , 0 , \fsl_2 , \fg_2 = \der \AA.
\end{array}$$

\medskip
This leaves three nilpotent orbits of $\fsl_6$ not in one of these series.
These correspond to the partitions $(51)$, $(321)$ and $(3111)$ and propagate
as follows: 

$$\begin{array}{lll}
\fg^2\fg_Q^2\qquad & [(51),(51,-),A_4]\qquad & \dim\cO_a = 12a+4, \\
 & & \dim  \fr (a) =3a, \\ & & \fh (a)= 0, 0, \fsl_3,\hspace{50mm} \\
 & & \\
\fg\fg_Q & [(321),(321,-),A_2+A_1]\hspace{15mm} & \dim\cO_a =9a+4, \\
 & & \dim  \fr (a) =5a+1, \\
 & & \fh (a)= 0, \fsl_2, \fsl_4, \\
 & & \\
\fg^2 & [(3111),(3111,-),A_2] & \dim\cO_a = 8a+2, \\
 & & \dim  \fr (a) =5a+1, \\
 & & \fh (a)= \fsl_3, \fso_6, \fsl_6. 
\end{array}$$

\subsection{The Severi series $\fg(\AA,\CC)$}
Here the Lie algebras $\fg$, and the number of positive roots $N$, 
parametrized by $a$ are:
\[ \begin{array}{ccccc}
a & 1 & 2 & 4 & 8 \\
\fg & \fsl_3 & 2\fsl_3 & \fsl_6 & \fe_6 \\
N & 3 & 6  & 15 & 36
\end{array} \]
The nilpotent orbits of $\fsl_3$ correspond to partitions of three and the nilpotent
orbits of $\fsl_6$ correspond to partitions of six. Given a partition
of three we construct a partition of six by repeating each part twice.
There are two dual preferred representations $V$ and $V^*$, of dimension $3a+3$, where
$V$ can be identified with the Jordan algebra $\cJ_3(\AA)$.
We obtain two series (we left a $?$ for the non-simple case, which has no standard label):

$$\begin{array}{lll}
V & [(21),(?),(2211),2A_2]\qquad &  \dim\cO_a = 4a, \\
 & & \dim\fr (a)=3a, \\
 & & \fh (a)= 0, 0, \fsl_2, \fg_2 = \der \AA, \\
 & & \\
\fg_Q=VV^* & [(3),(?),(33),2A_1] & \dim\cO_a =6a, \\
 & & \dim\fr (a)=2a, \\ 
 & & \fh (a)=0, \CC, \fsl_2, \fso_7=\fl(\AA).
\end{array}$$
 
\subsection{The sub-Severi series $\fg(\AA,\RR)$}
This is the series $\fg(\AA,\RR)=\fsl_2, \fsl_3, \fsp_6, \ff_4$, with its 
preferred representation $W=\cJ_3(\AA)_0$ of dimension $3a+2$: the space of traceless
matrices in $\cJ_3(\AA)$. This leads to the following series of orbits:
 
$$\begin{array}{lll}
\fg_Q=W\qquad & [(2),(3),(3,-),\tilde A_2]\qquad  & \dim\cO_a =4a-2, \\
 & & \dim  \fr (a) =a, \\
 & & \fh (a) = 0, 0, \fsl_2 , \fg_2 =\der \AA.
\end{array}$$

 \bigskip

\section{Beyond the exceptional Lie algebras}

\subsection{General dimension formulas}
There are four nonzero nilpotent orbits occuring in all simple Lie algebras
of rank greater than two (and also $\fg_2$):
\begin{description}
\item[1] the {\sl regular nilpotent orbit}, which is the open orbit in the nilpotent cone, 
\item[2] the {\sl subregular nilpotent orbit}, which is the open orbit in the boundary 
of the regular orbit,
\item[3] the {\sl minimal nilpotent orbit}, which we call $\cO_{ad}\subset \fg$
(we often work with its projectivization   $X_{ad}\subset\BP \fg$),
and in this paper is denoted simply $\fg$ as the marked Dynkin diagram
corresponds to the adjoint representation,
\item[4] the orbit  whose projectivization we called $\s_{(1)}(X_{ad})$ in
\cite{LMmagic}. 
\end{description}
Panyushev, in \cite{pan}, calls this last orbit $\BO$, but because of our 
usage of $\BO$ to denote the octonions, we will denote it by $\cO_{\s_{(1)}(X_{ad})}$
or $\fg_2$, since its marked Dynkin diagram gives the weight of $\fg_2$.  
Note that Panyushev only observes this orbit when the adjoint representation is fundamental, 
where it corresponds to the diagram marked with a $1$ over nodes adjacent to the node 
of the adjoint representation, and zeros elsewhere. Geometrically $\cO_{\s_{(1)}(X_{ad})}$ may
be described as either the union of tangent lines to the contact
distribution on $X_{ad}$, or as the 
closure of the set of points in $\BP \fg$ lying
on a two-parameter family of secant lines, see \cite{LMmagic}.

\smallskip
The dimension of the regular nilpotent orbit has a simple expression,
either the number of roots, or the dimension of $\fg$ minus the rank
of $\fg$, and the subregular orbit, being of codimension two in
the closure of the regular orbit, inherits a dimension formula.

Nevertheless, when we study orbits in series, we see that from the series perspective,
the properties of being regular and subregular are not good ones. What
happens instead is that the regular and subregular orbits of the fixed
algebra in a series gives rise to a series of orbits which
in general are not regular or subregular. 
This is not surprising as the dimension of the regular and subregular
orbits grow like the square of the parameter parametrizing the algebras
(as do  the dimensions of the algebras themselves), while we insist that
the nilpotent orbits in series have linear dimension formulas.

The starting point of Vogel's conjectured {\sl universal Lie 
algebra} was an attempt to construct a category with analogs of the Casimir, the bracket, 
the Killing form and  the Jacobi identity, that dominates the category of modules of any 
simple Lie algebra. It leads to a parametrization of the simple Lie algebras by a projective
plane, whose barycentric coordinate is the eigenvalue of the Casimir operator on the
adjoint representation, 
and the scaling is by the length of the longest root. (See \cite{vog2,del} for these 
parameters, and \cite{LMtrial,LMpop} for the relation with the triality model.)

Remarkably, the minimal nilpotent orbit has a nice 
dimension formula in the spirit of Vogel's work. This was first observed by W. Wang
(\cite{wang}, independently of this interpretation). 

\begin{proposition} Let $\fg$ be a complex simple Lie algebra. After an invariant 
quadratic form has been chosen, let $\sqrt a$ denote the length of the longest root,  and 
$C$ the Casimir eigenvalue for $\fg$. Then
$$\dim\cO_{ad} = \frac{2C}a -2. $$
\end{proposition}
 
Wang's formula is actually $\dim\cO_{ad} =2\check{h}-2$, where $\check{h}$
denotes the dual Coxeter number. But, once we have fixed an invariant 
scalar product on the root lattice, we can write
$$\frac{C}{a}=\frac{\langle\tilde{\a}+2\rho,\tilde{\a}
\rangle}{\langle\tilde{\a},\tilde{\a}\rangle}=1+2
\frac{\langle\rho,\tilde{\a}\rangle}{\langle\tilde{\a},\tilde{\a}
\rangle}=\check{h},$$
the last equality being a definition. Here $\tilde{a}$ denotes the
highest root, and $2\rho$ the sum of all positive roots. 
Note that this is just a linear formula, while Vogel's dimension formulas for the modules 
are much more complicated. 

\smallskip
In \cite{LMmagic} we discuss two other series of nilpotent orbits that
are not completely general. We revert to the notation of \cite{LMmagic}, discussing
the projectivizations of the orbit closures in $\BP \fg$. In particular
the dimension of the corresponding orbit closure is one more than that of its
projectivization.
\begin{description}
\item[5]  $\cO_{\s_{(3)}(X_{ad})}$, : this orbit occurs in the
exceptional series (with label $\fg^2$) and the $\fsl$ series. 
Geometrically, it is the union of tangent lines to $X_{ad}$ that are
tangent to the quartic cone inside each hyperplane in the contact distribution.
In terms of weighted Dynkin diagrams, one marks the adjoint
nodes with a $2$ and puts zeros elsewhere.
\item[6] $\cO_{\s_{Q}(X_{ad})}$ : here $Q$ denotes an unextendable quadric on $X_{ad}$. This series 
occurs in the exceptional series (with label $\fg_Q$)
and there are two different series of such orbits in the $\fso$-series (even three
for $\fso_8$, but they are all isomorphic). In each series
the dimension of $Q$ is a linear function of the parameter parametrizing the
series. Geometrically, these orbits  are obtained by taking a uniruling of $X_{ad}$
by unextendable quadrics and taking the union of their projective spans.
In terms of weighted Dynkin diagrams, one marks with a $1$ the node such that, when erased, 
the connected component of the node marked for the adjoint representation is a marked Dynkin 
diagram corresponding to a quadric hypersurface, see \cite{LMmagic}.
\end{description}

The dimensions of these orbits for the exceptional series were computed in \cite{LMmagic}.
For the classical series they can be extracted from \cite{carter}, and we get : 

\begin{corollary} 
$$\dim \cO_{\s_{(1)}(X_{ad})} =\frac{4C}a -5, \quad \dim \cO_{\s_{(3)}(X_{ad})}
=\frac{4C}a -9, \quad \dim \cO_{\s_{Q}(X_{ad})} =  \frac{4C}a-\dim Q -5.$$
\end{corollary}

\medskip\noindent {\sl Remark}.
Note that  in the classical series, if one extends a partition by
zero, one obtains a series of nilpotent orbits in our sense, in that
the dimensions of the orbits are given as linear functions of the
parameters. More precisely, we have:

\begin{proposition} Fix $f$, and respectively let $\fg_f=
 \fsl_f,\fso_f,\fsp_{2f}$ and $\fg(t)=
 \fg_{f+t}$. Let $\cO$ be a nilpotent orbit in $\fg_f$.
Let $r_i$ denote the number of elementary divisors with exponent $i$ 
in the partition defining $\cO$ (following \cite{carter}). Let
$\cO_t\subset \fg_t$ be the corresponding orbit with $r_1(t)= t+r_1$
and all the other $r_i$'s the same.
Then $\dim \cO_t$ is a linear function of $t$. More precisely we have
$$\begin{array}{llcll}
(1)\qquad & \dim \cO_t & =  & 2t(f-(r_1+\cdots + r_n)) +\dim \cO &\qquad {\rm sl-case} \\
(2) & \dim \cO_t & =  & t(f-(r_1+\cdots + r_n)-\frac 34) +\dim \cO &  \qquad{\rm so-case} \\
(3) & \dim \cO_t & =  & 2t(f-(r_1+\cdots + r_n)+\frac 34) +\dim \cO & \qquad {\rm sp-case}
\end{array}$$ 
\end{proposition}

\medskip\noindent
As with the exceptional series, these orbits also share
a common geometry. Their desingularizations by vector bundles
$E\ra G/P$ are such that the spaces $G/P$ have uniform geometric 
interpretations,
which are obvious in the classical cases, and can be understood
uniformly in terms of their shadows on the adjoint varieties $X_{ad}\subset\BP \fg$.
The $P$-modules defining $E$ also have uniform interpretations in
terms of Tits geometries.
\medskip
 
\subsection{The generalized magic square}
This is the following $3\times 3$ square, with parameters $n\ge 4$ and $a,b=1,2,4$:
$$\begin{array}{rccc}
 & a=1 & a=2 & a=4 \\
b=1 & \fso_n & \fsl_n & \fsp_{2n} \\
b=2 & \fsl_n & 2\fsl_n & \fsl_{2n} \\
b=4 & \fsp_{2n} & \fsl_{2n} & \fso_{4n}
\end{array}$$

\medskip
Recall from \cite{carter} or \cite{colmc} \S 5.1,  that nilpotent orbits in $\fsl_n$,
respectively $\fsp_{2n},\fso_{2n+1},\fso_{2n}$ are in  one to one
correspondence with partitions $(d_1,\ldots, d_n)$ of $n$,  
respectively partitions of $2n$ in which odd parts occur with
even multiplicity, partitions of $2n+1$ in which even parts occur
with even multiplicity, partitions of $2n$ in which even parts occur
with even multiplicity (with a slight modification for partitions with
only even parts). We let $r_i$ be the number of times $i$ occurs in the partition.

\begin{proposition}
 For every nilpotent orbit $\cO_{1,1}$ in $\fso_n$, there is a  nilpotent orbit 
$\cO_{a,b}$ for each element of the generalized magic square,  whose dimension is a 
bilinear function
of $a$ and $b$.  More precisely, let $(r_1\hd r_n)$ be as above for the partition 
parametrizing the orbit $\cO_{1,1}\subset \fso_n$. Then
$$
\dim \cO_{a,b}= \frac{ab}2 (n^2- \Sigma_i(\Sigma_{j\ge i}r_j)^2 -n+\Sigma_{i\ odd}r_i)
+(a+b-2) (n-\Sigma_{i\ odd}r_i). 
$$
\end{proposition}

More generally one can take a nilpotent orbit for any algebra in the square and 
extend it across and below to get a bilinear function in $a,b$. 

\begin{proof}
Given a partition of $n$ admissible for $\fso_n$, just use it
as a partition for $\fsl_n$ to get a nilpotent orbit. Given a partition
of $\fsl_n$, double it to get an admissible partition for $\fsp_{2n}$.
Given a partition for $\fsl_n$ use it twice to get a partition for
$2\fsl_n$. Given two partitions of length $n$ (parametrizing 
a nilpotent orbit in $2\fsl_n$) put them together to get a partition
for $\fsl_{2n}$. Given a partition admissible for $\fsp_{2n}$, just use
it to get a partition for $\fsl_{2n}$. Given a partition of $\fsl_{2n}$,
double it to get an admissible partition for $\fso_{4n}$.
These are the partitions we use to define $\cO_{a,b}$ from $\cO_{1,1}$. Note that
the process is symmetric in how one moves across the chart.

Now the proof just consists in checking that for each value of $a,b$, the dimension 
given by the above formula is consistent with those given in \cite{colmc} for the classical
Lie algebras. \end{proof}

Note that this works also for the $n=3$ chart including the exceptional
groups, and also that one can begin anywhere in the chart to get orbits
to the right and below. Finally, specializing to each row, one gets
{\it linear} functions of $a$ for the dimensions.

\medskip\noindent {\sl Examples}.
 
1. The regular nilpotent orbit in $\fso_n$ induces a series with
$$\dim \cO_{a,b}= \frac{ab}2(n^2-n-1+\epsilon) + (a+b+2)(n-\epsilon),
$$
where $\epsilon =1$ if $n$ is odd and $0$ if $n$ is even.

2. (A magical orbit.) Consider the partition $(3,1\hd 1)$ and the
resulting three parameter family of orbits $\cO_{a,b,n}$. We get the following 
dimension formula, which has the very nice property of being linear in each of the 
parameters:
$$\dim\cO_{a,b,n}= 2(ab(n-2) +  a+b-2).$$
Note that this orbit is {\it universal} in that it occurs in all simple Lie algebras:
in the exceptional and sub-exceptional cases this is the series labeled $\fg_Q^2$, in the
Severi case the series labeled $\fg_Q=VV^*$, and in the sub-Severi case the series 
labeled $\fg_Q=W$.

\vspace{1cm}

{\small
\noindent {\sc Joseph M. Landsberg}, School of Mathematics, 
Georgia Institute of Technology, Atlanta, GA 30332-0160, USA

\noindent {\rm E-mail}: jml@math.gatech.edu 

\medskip\noindent
{\sc Laurent Manivel},  Institut Fourier, UMR 5582 du CNRS 
Universit\'e Grenoble I, BP 74, 38402 Saint Martin d'H\`eres cedex, FRANCE

\noindent {\rm E-mail}: Laurent.Manivel@ujf-grenoble.fr 

\medskip\noindent {\sc Bruce W. Westbury}, Mathematics Institute, University of Warwick, 
Coventry CV4 7AL, UK 

\noindent {\rm E-mail}: bww@maths.warwick.ac.uk  }

\end{document}